\documentclass[final]{siamltex}
\usepackage{fancybox,varioref}
\usepackage{amsmath,amssymb}
\usepackage[dvips]{graphicx}
\newcommand{\ul}{\underline}
\newcommand{\ol}{\overline}
\newcommand{\ts}{\textstyle}

\newcommand{\Real}{\ensuremath{\mathrm{R\mspace{-1mu}e}}}
\newcommand{\Imag}{\ensuremath{\mathrm{I\mspace{-2mu}m}}}
\newcommand{\hc}{\ensuremath{\hat{c}}}
\newcommand{\hck}{\ensuremath{\hc^{(k)}}}
\newcommand{\hcko}{\ensuremath{\hc^{(k+1)}}}
\newcommand{\hcs}{\ensuremath{\hc^{(0)}}}
\newcommand{\czrw}{\ensuremath{c_\mathrm{rw}^\star}}

\newcommand{\aR}{\ensuremath{\alpha_{\scriptscriptstyle R}}}
\newcommand{\aI}{\ensuremath{\alpha_{\scriptscriptstyle I}}}
\newcommand{\aRl}{\ensuremath{\aR(\lam)}}
\newcommand{\aIl}{\ensuremath{\aI(\lam)}}
\newcommand{\daRz}{\ensuremath{\dot{\alpha}_{\scriptscriptstyle R}^{\star}}}

\newcommand{\daRlz}{\ensuremath{\dot{\alpha}_{\scriptscriptstyle R}(\lamz)}}
\newcommand{\daIlz}{\ensuremath{\dot{\alpha}_{\scriptscriptstyle I}(\lamz)}}
\newcommand{\beR}{\ensuremath{\beta_{\scriptscriptstyle\mathrm{R}}}}
\newcommand{\beI}{\ensuremath{\beta_{\scriptscriptstyle\mathrm{I}}}}
\newcommand{\bRl}{\ensuremath{\beR(\lam)}}
\newcommand{\bIl}{\ensuremath{\beI(\lam)}}
\newcommand{\bRlz}{\ensuremath{\beR(\lamz)}}
\newcommand{\bIlz}{\ensuremath{\beI(\lamz)}}
\newcommand{\kapl}{\ensuremath{\kappa(\lam)}}
\newcommand{\kaplz}{\ensuremath{\kappa(\lamz)}}

\newcommand{\ep}{\ensuremath{\varepsilon}}
\newcommand{\rep}{\ensuremath{\tfrac{1}{\ep}}}
\newcommand{\rept}{\ensuremath{\tfrac{1}{\ep^2}}}
\newcommand{\bgam}{\ensuremath{\boldsymbol{\gamma}}}

\newcommand{\bsig}{\ensuremath{\boldsymbol{\sigma}}}
\newcommand{\bsigrw}{\ensuremath{\bsig^\mathrm{rw}}}
\newcommand{\bsigpns}{\ensuremath{\bsig_{\scriptscriptstyle +}^\mathrm{ns}}}
\newcommand{\bsigmns}{\ensuremath{\bsig_{\scriptscriptstyle -}^\mathrm{ns}}}
\newcommand{\bsigRns}{\ensuremath{\bsig_{\scriptscriptstyle\mathrm{R}}^\mathrm{ns}}}
\newcommand{\bsigIns}{\ensuremath{\bsig_{\scriptscriptstyle\mathrm{I}}^\mathrm{ns}}}
\newcommand{\tl}{\ensuremath{(\theta,\lam)}}
\newcommand{\tlz}{\ensuremath{(\theta,\lamz)}}
\newcommand{\tp}{\ensuremath{(\theta,\phi)}}
\newcommand{\tep}{\ensuremath{(\theta;\ep)}}
\newcommand{\tpe}{\ensuremath{(\theta,\phi;\ep)}}
\newcommand{\etas}{\ensuremath{\eta^{(0)}}}
\newcommand{\etak}{\ensuremath{\eta^{(k)}}}
\newcommand{\etako}{\ensuremath{\eta^{(k+1)}}}
\newcommand{\etp}{\ensuremath{\eta\tp}}
\newcommand{\zetp}{\ensuremath{\zeta\tp}}
\newcommand{\half}{\ensuremath{\frac12}}
\newcommand{\thalf}{\ensuremath{{\ts\frac12}}}
\newcommand{\np}{\ensuremath{n_+}}
\newcommand{\nm}{\ensuremath{n_-}}
\newcommand{\npm}{\ensuremath{n_{\pm}}}
\newcommand{\mA}{\ensuremath{\mathsf{A}}}
\newcommand{\mAt}{\ensuremath{\mA(\theta)}}

\newcommand{\mBh}{\ensuremath{\widehat{\mathsf{B}}}}
\newcommand{\mBhl}{\ensuremath{\mBh(\lam)}}

\newcommand{\mD}{\ensuremath{\mathsf{D}}}
\newcommand{\mDL}{\ensuremath{\mD_L}}
\newcommand{\mE}{\ensuremath{\mathsf{E}}}
\newcommand{\mEl}{\ensuremath{\mE(\lam)}}

\newcommand{\dElz}{\ensuremath{\dot{\mE}(\lamz)}}
\newcommand{\mEm}{\ensuremath{\mE_-}}
\newcommand{\mEp}{\ensuremath{\mE_+}}
\newcommand{\mEpm}{\ensuremath{\mE_{\pm}}}
\newcommand{\mEz}{\ensuremath{\mE^\star}}

\newcommand{\wE}{\ensuremath{\widetilde{\mE}}}
\newcommand{\wEz}{\ensuremath{\wE^\star}}
\newcommand{\wEl}{\ensuremath{\wE(\lam)}}

\newcommand{\mEh}{\ensuremath{\widehat{\mE}}}
\newcommand{\mEhz}{\ensuremath{\mEh^\star}}
\newcommand{\mEhl}{\ensuremath{\mEh(\lam)}}
\newcommand{\mH}{\ensuremath{\mathsf{H}}}

\newcommand{\mI}{\ensuremath{\mathsf{I}}}

\newcommand{\mJ}{\ensuremath{\mathsf{J}}}
\newcommand{\mJRss}{\ensuremath{\mJ_{\scriptscriptstyle\mathrm{R}}^\mathrm{ss}}}
\newcommand{\mJIss}{\ensuremath{\mJ_{\scriptscriptstyle\mathrm{I}}^\mathrm{ss}}}
\newcommand{\mJrw}{\ensuremath{\mJ^\mathrm{rw}}}
\newcommand{\mJRRrw}{\ensuremath{\mJrw_{\scriptscriptstyle\mathrm{RR}}}}
\newcommand{\mJRIrw}{\ensuremath{\mJrw_{\scriptscriptstyle\mathrm{RI}}}}
\newcommand{\mJIRrw}{\ensuremath{\mJrw_{\scriptscriptstyle\mathrm{IR}}}}
\newcommand{\mJIIrw}{\ensuremath{\mJrw_{\scriptscriptstyle\mathrm{II}}}}
\newcommand{\mJCrw}{\ensuremath{\mJrw_{\scriptscriptstyle\C}}}
\newcommand{\mJpprw}{\ensuremath{\mJrw_{\scriptscriptstyle ++}}}
\newcommand{\mJpmrw}{\ensuremath{\mJrw_{\scriptscriptstyle +-}}}
\newcommand{\mJmprw}{\ensuremath{\mJrw_{\scriptscriptstyle -+}}}
\newcommand{\mJmmrw}{\ensuremath{\mJrw_{\scriptscriptstyle --}}}
\newcommand{\mP}{\ensuremath{\mathsf{P}}}
\newcommand{\mPt}{\ensuremath{\mP(\theta)}}
\newcommand{\mPl}{\ensuremath{\mP(\lam)}}
\newcommand{\mPtl}{\ensuremath{\mP(\theta;\lam)}}

\newcommand{\mPz}{\ensuremath{\mP^\star}}
\newcommand{\mPzt}{\ensuremath{\mPz(\theta)}}
\newcommand{\dP}{\ensuremath{\dot{\mP}}}
\newcommand{\dPt}{\ensuremath{\dP(\theta)}}
\newcommand{\dPtl}{\ensuremath{\dP(\theta;\lam)}}
\newcommand{\mPp}{\ensuremath{\mP_+}}
\newcommand{\mPpt}{\ensuremath{\mPp(\theta)}}
\newcommand{\mPm}{\ensuremath{\mP_-}}
\newcommand{\mPmt}{\ensuremath{\mPm(\theta)}}
\newcommand{\mPpm}{\ensuremath{\mP_{\pm}}}
\newcommand{\mPpmtpt}{\ensuremath{\mPpm(\theta+\tpi)}}

\newcommand{\dPzt}{\ensuremath{\dP^\star(\theta)}}

\newcommand{\mPpmt}{\ensuremath{\mPpm(\theta)}}

\newcommand{\mS}{\ensuremath{\mathsf{S}}}
\newcommand{\mSM}{\ensuremath{\mS_{\scriptscriptstyle M}}}
\newcommand{\mSLM}{\ensuremath{\mS_{\scriptscriptstyle L,M}}}

\newcommand{\mT}{\ensuremath{\mathsf{T}}}
\newcommand{\mTa}{\ensuremath{\mT_{\scriptscriptstyle\mathrm{a}}}}
\newcommand{\mTs}{\ensuremath{\mT_{\scriptscriptstyle\mathrm{s}}}}

\newcommand{\ZZ}{\ensuremath{\mathsf{O}}}
\newcommand{\C}{\ensuremath{\mathbb{C}}}
\newcommand{\CL}{\ensuremath{\C^L}}

\newcommand{\So}{\ensuremath{\mathbb{S}^1}}
\newcommand{\Z}{\ensuremath{\mathbb{Z}}}
\newcommand{\R}{\ensuremath{\mathbb{R}}}
\newcommand{\RL}{\ensuremath{\R^L}}
\newcommand{\Rt}{\ensuremath{\R^2}}
\newcommand{\Rn}{\ensuremath{\R^n}}

\newcommand{\Rnt}{\ensuremath{\R^{n-2}}}
\newcommand{\Rnth}{\ensuremath{\R^{n-3}}}

\newcommand{\Rnn}{\ensuremath{\R^{n\times n}}}
\newcommand{\Rtt}{\ensuremath{\R^{2\times 2}}}

\newcommand{\Rntnt}{\ensuremath{\R^{(n-2)\times(n-2)}}}

\newcommand{\Rnthnth}{\ensuremath{\R^{(n-3)\times(n-3)}}}

\newcommand{\Rnpnp}{\ensuremath{\R^{\np\times\np}}}

\newcommand{\Rnmnm}{\ensuremath{\R^{\nm\times\nm}}}
\newcommand{\Rnnpm}{\ensuremath{\R^{n\times\npm}}}
\newcommand{\bF}{\ensuremath{\boldsymbol{F}}}
\newcommand{\bFxl}{\ensuremath{\bF(\bx,\lam)}}
\newcommand{\bFxtl}{\ensuremath{\bF(\bx,t,\lam)}}
\newcommand{\bFxzlz}{\ensuremath{\bF(\bxz,\lamz)}}
\newcommand{\bFxzll}{\ensuremath{\bF(\bxzl,\lam)}}
\newcommand{\bFvtl}{\ensuremath{\bF(\bv,\theta,\lam)}}
\newcommand{\bFvztlz}{\ensuremath{\bF(\bvzt,\lamz)}}
\newcommand{\bFvztll}{\ensuremath{\bF(\bvztl,\lam)}}
\newcommand{\bFvzttlz}{\ensuremath{\bF(\bvzt,\theta,\lamz)}}
\newcommand{\bFvztltl}{\ensuremath{\bF(\bvztl,\theta,\lam)}}

\newcommand{\bG}{\ensuremath{\boldsymbol{G}}}
\newcommand{\bGh}{\ensuremath{\boldsymbol{\widehat{G}}}}
\newcommand{\bGhd}{\ensuremath{\bGh^{\dag}}}
\newcommand{\bGhdd}{\ensuremath{\bGh^{\ddag}}}
\newcommand{\wbG}{\ensuremath{\boldsymbol{\widetilde{G}}}}
\newcommand{\wbGd}{\ensuremath{\wbG^{\dag}}}
\newcommand{\wbGdd}{\ensuremath{\wbG^{\ddag}}}
\newcommand{\Go}{\ensuremath{\mathring{G}}}
\newcommand{\God}{\ensuremath{\Go^{\dag}}}
\newcommand{\Godd}{\ensuremath{\Go^{\ddag}}}

\newcommand{\bKh}{\ensuremath{\boldsymbol{\widehat{K}}}}
\newcommand{\bKhd}{\ensuremath{\bKh^{\dag}}}
\newcommand{\wbK}{\ensuremath{\boldsymbol{\widetilde{K}}}}
\newcommand{\wbKd}{\ensuremath{\wbK^{\dag}}}
\newcommand{\Ko}{\ensuremath{\mathring{K}}}
\newcommand{\Kod}{\ensuremath{\Ko^{\dag}}}
\newcommand{\bQ}{\ensuremath{\boldsymbol{Q}}}
\newcommand{\bQC}{\ensuremath{\bQ_{\scriptscriptstyle\C}}}
\newcommand{\QC}{\ensuremath{Q_{\scriptscriptstyle\C}}}
\newcommand{\bQI}{\ensuremath{\bQ_{\scriptscriptstyle\mathrm{I}}}}
\newcommand{\QI}{\ensuremath{Q_{\scriptscriptstyle\mathrm{I}}}}

\newcommand{\oa}{\ensuremath{\mathring{a}}}
\newcommand{\ba}{\ensuremath{\boldsymbol{a}}}
\newcommand{\baz}{\ensuremath{\ba^\star}}
\newcommand{\bazt}{\ensuremath{\baz(\theta)}}
\newcommand{\bazp}{\ensuremath{\baz(\phi)}}
\newcommand{\bah}{\ensuremath{\boldsymbol{\widehat{a}}}}
\newcommand{\bahz}{\ensuremath{\bah^\star}}
\newcommand{\bahzt}{\ensuremath{\bahz(\theta)}}

\newcommand{\bahzp}{\ensuremath{\bahz(\phi)}}
\newcommand{\badhzp}{\ensuremath{\boldsymbol{\widehat{\dot{a}}}^\star(\phi)}}
\newcommand{\dba}{\ensuremath{\boldsymbol{\dot{a}}}}
\newcommand{\dbaz}{\ensuremath{\dba^\star}}
\newcommand{\dbazt}{\ensuremath{\dbaz(\theta)}}
\newcommand{\dbazp}{\ensuremath{\dbaz(\phi)}}
\newcommand{\wba}{\ensuremath{\boldsymbol{\widetilde{a}}}}
\newcommand{\ob}{\ensuremath{\mathring{b}}}
\newcommand{\bb}{\ensuremath{\boldsymbol{b}}}
\newcommand{\bbz}{\ensuremath{\bb^\star}}
\newcommand{\bbh}{\ensuremath{\boldsymbol{\widehat{b}}}}
\newcommand{\wbb}{\ensuremath{\boldsymbol{\widetilde{b}}}}

\newcommand{\bd}{\ensuremath{\boldsymbol{d}}}
\newcommand{\be}{\ensuremath{\boldsymbol{e}}}
\newcommand{\bens}{\ensuremath{\be^\mathrm{ns}}}
\newcommand{\beRns}{\ensuremath{\bens_{\scriptscriptstyle\mathrm{R}}}}
\newcommand{\beIns}{\ensuremath{\bens_{\scriptscriptstyle\mathrm{I}}}}
\newcommand{\bef}{\ensuremath{\be_1}}
\newcommand{\bet}{\ensuremath{\be_2}}
\newcommand{\beh}{\ensuremath{\boldsymbol{\widehat{e}}}}

\newcommand{\ben}{\ensuremath{\be_n}}
\newcommand{\berw}{\ensuremath{\be^\mathrm{rw}}}
\newcommand{\beRrw}{\ensuremath{\berw_{\scriptscriptstyle\mathrm{R}}}}
\newcommand{\beIrw}{\ensuremath{\berw_{\scriptscriptstyle\mathrm{I}}}}
\newcommand{\bff}{\ensuremath{\boldsymbol{f}}}

\newcommand{\bft}{\ensuremath{\bff(\theta)}}

\newcommand{\bg}{\ensuremath{\boldsymbol{g}}}
\newcommand{\bgt}{\ensuremath{\bg(\theta)}}
\newcommand{\bgm}{\ensuremath{\bg_-}}
\newcommand{\bgmt}{\ensuremath{\bgm(\theta)}}
\newcommand{\bgp}{\ensuremath{\bg_+}}
\newcommand{\bgpt}{\ensuremath{\bgp(\theta)}}

\newcommand{\ho}{\ensuremath{\mathring{h}}}

\newcommand{\bhh}{\ensuremath{\boldsymbol{\widehat{h}}}}
\newcommand{\wbh}{\ensuremath{\boldsymbol{\widetilde{h}}}}

\newcommand{\bq}{\ensuremath{\boldsymbol{q}}}
\newcommand{\bqR}{\ensuremath{\bq_{\scriptscriptstyle \mathrm{R}}}}
\newcommand{\bqI}{\ensuremath{\bq_{\scriptscriptstyle \mathrm{I}}}}
\newcommand{\bqRk}{\ensuremath{\bqR^{(k)}}}
\newcommand{\bqIk}{\ensuremath{\bqI^{(k)}}}
\newcommand{\bqp}{\ensuremath{\bq_{\scriptscriptstyle +}}}
\newcommand{\bqm}{\ensuremath{\bq_{\scriptscriptstyle -}}}
\newcommand{\bqpk}{\ensuremath{\bqp^{(k)}}}
\newcommand{\bqmk}{\ensuremath{\bqm^{(k)}}}
\newcommand{\br}{\ensuremath{\boldsymbol{r}}}
\newcommand{\brR}{\ensuremath{\br_{\scriptscriptstyle\mathrm{R}}}}
\newcommand{\brRk}{\ensuremath{\brR^{(k)}}}
\newcommand{\brI}{\ensuremath{\br_{\scriptscriptstyle\mathrm{I}}}}
\newcommand{\brIk}{\ensuremath{\brI^{(k)}}}
\newcommand{\brp}{\ensuremath{\br_{\scriptscriptstyle +}}}
\newcommand{\brm}{\ensuremath{\br_{\scriptscriptstyle -}}}
\newcommand{\brpk}{\ensuremath{\brp^{(k)}}}
\newcommand{\brmk}{\ensuremath{\brm^{(k)}}}

\newcommand{\brk}{\ensuremath{\br^{(k)}}}

\newcommand{\wbrktp}{\ensuremath{\widetilde{\br}^{(k)}\tp}}
\newcommand{\hrktp}{\ensuremath{\widehat{r}^{(k)}\tp}}
\newcommand{\bs}{\ensuremath{\boldsymbol{s}}}
\newcommand{\hs}{\ensuremath{\hat{\bs}}}
\newcommand{\hsk}{\ensuremath{\hs^{(k)}}}
\newcommand{\hsf}{\ensuremath{\hs^{(0)}}}

\newcommand{\bu}{\ensuremath{\boldsymbol{u}}}
\newcommand{\buz}{\ensuremath{\bu^\star}}

\newcommand{\bw}{\ensuremath{\boldsymbol{w}}}
\newcommand{\bwt}{\ensuremath{\bw(\theta)}}

\newcommand{\bwm}{\ensuremath{\bw_-}}
\newcommand{\bwmt}{\ensuremath{\bwm(\theta)}}
\newcommand{\bwp}{\ensuremath{\bw_+}}

\newcommand{\bwpt}{\ensuremath{\bwp(\theta)}}

\newcommand{\bv}{\ensuremath{\boldsymbol{v}}}
\newcommand{\bvR}{\ensuremath{\bv_{\scriptscriptstyle \mathrm{R}}}}
\newcommand{\bvRk}{\ensuremath{\bvR^{(k)}}}
\newcommand{\bvRko}{\ensuremath{\bvR^{(k+1)}}}
\newcommand{\bvRs}{\ensuremath{\bvR^{(0)}}}
\newcommand{\bvI}{\ensuremath{\bv_{\scriptscriptstyle \mathrm{I}}}}
\newcommand{\bvIk}{\ensuremath{\bvI^{(k)}}}
\newcommand{\bvIko}{\ensuremath{\bvI^{(k+1)}}}
\newcommand{\bvIs}{\ensuremath{\bvI^{(0)}}}

\newcommand{\bvt}{\ensuremath{\bv(\theta)}}
\newcommand{\bvz}{\ensuremath{\bv^\star}}
\newcommand{\bvzt}{\ensuremath{\bvz(\theta)}}
\newcommand{\bvztl}{\ensuremath{\bvz(\theta;\lam)}}

\newcommand{\bx}{\ensuremath{\boldsymbol{x}}}
\newcommand{\bxz}{\ensuremath{\bx^\star}}

\newcommand{\hy}{\ensuremath{\widehat{y}}}
\newcommand{\by}{\ensuremath{\boldsymbol{y}}}

\newcommand{\byh}{\ensuremath{\boldsymbol{\widehat{y}}}}

\newcommand{\byht}{\ensuremath{\byh(\theta)}}
\newcommand{\byhtp}{\ensuremath{\byh\tp}}
\newcommand{\wty}{\ensuremath{\widetilde{\by}}}
\newcommand{\wtyLM}{\ensuremath{\wty_{\scriptscriptstyle L,M}}}
\newcommand{\wtyz}{\ensuremath{\wty^{\star}}}
\newcommand{\wtyt}{\ensuremath{\wty(\theta)}}
\newcommand{\wtytp}{\ensuremath{\wty\tp}}

\newcommand{\yo}{\ensuremath{\mathring{y}}}
\newcommand{\yotp}{\ensuremath{\yo\tp}}
\newcommand{\bz}{\ensuremath{\boldsymbol{z}}}
\newcommand{\bzs}{\ensuremath{\bz^{(0)}}}
\newcommand{\bzR}{\ensuremath{\bz_{\scriptscriptstyle \mathrm{R}}}}
\newcommand{\bzRk}{\ensuremath{\bzR^{(k)}}}
\newcommand{\bzRko}{\ensuremath{\bzR^{(k+1)}}}
\newcommand{\bzI}{\ensuremath{\bz_{\scriptscriptstyle \mathrm{I}}}}
\newcommand{\bzIk}{\ensuremath{\bzI^{(k)}}}
\newcommand{\bzIko}{\ensuremath{\bzI^{(k+1)}}}
\newcommand{\bzp}{\ensuremath{\bz_{\scriptscriptstyle +}}}
\newcommand{\bzm}{\ensuremath{\bz_{\scriptscriptstyle -}}}
\newcommand{\bzps}{\ensuremath{\bzp^{(0)}}}
\newcommand{\bzms}{\ensuremath{\bzm^{(0)}}}
\newcommand{\bzpk}{\ensuremath{\bzp^{(k)}}}
\newcommand{\bzmk}{\ensuremath{\bzm^{(k)}}}
\newcommand{\bzpko}{\ensuremath{\bzp^{(k+1)}}}
\newcommand{\bzmko}{\ensuremath{\bzm^{(k+1)}}}
\newcommand{\bzM}{\ensuremath{\bz_{\scriptscriptstyle M}}}
\newcommand{\bzz}{\ensuremath{\bz^\star}}
\newcommand{\bzh}{\ensuremath{\boldsymbol{\widehat{z}}}}
\newcommand{\bzhz}{\ensuremath{\boldsymbol{\widehat{z}}^\star}}
\newcommand{\wtz}{\ensuremath{\widetilde{\bz}}}
\newcommand{\wtzz}{\ensuremath{\widetilde{\bz}^\star}}
\newcommand{\wtzt}{\ensuremath{\wtz(\theta)}}

\newcommand{\zo}{\ensuremath{\mathring{z}}}

\newcommand{\bzt}{\ensuremath{\bz(\theta)}}
\newcommand{\bzht}{\ensuremath{\bzh(\theta)}}
\newcommand{\bztp}{\ensuremath{\bz\tp}}
\newcommand{\bzhtp}{\ensuremath{\bzh\tp}}
\newcommand{\bzst}{\ensuremath{\bzs(\theta)}}

\newcommand{\bzMt}{\ensuremath{\bzM(\theta)}}

\newcommand{\zz}{\ensuremath{\boldsymbol{0}}}
\newcommand{\ee}{\ensuremath{\mathrm{e}}}
\newcommand{\ii}{\ensuremath{\mathrm{i}}}
\newcommand{\dif}{\ensuremath{\mathrm{d}}}
\newcommand{\dt}{\ensuremath{\,\dif\theta}}
\newcommand{\dphi}{\ensuremath{\,\dif\phi}}
\newcommand{\ddt}{\ensuremath{\frac{\dif\;{}}{\dif\theta}}}
\newcommand{\dvdt}{\ensuremath{\frac{\dif\bv}{\dif\theta}}}

\newcommand{\dxdt}{\ensuremath{\frac{\dif\bx}{\dif t}}}
\newcommand{\pdt}{\ensuremath{\frac{\partial\;{}}{\partial\theta}}}
\newcommand{\pdp}{\ensuremath{\frac{\partial\;{}}{\partial\phi}}}
\newcommand{\pbzdttp}{\ensuremath{\frac{\partial\bz}{\partial\theta}\tp}}
\newcommand{\pbzdptp}{\ensuremath{\frac{\partial\bz}{\partial\phi}\tp}}

\newcommand{\xt}{\ensuremath{(x,t)}}
\newcommand{\pdx}{\ensuremath{\frac{\partial\;{}}{\partial x}}}
\newcommand{\pudt}{\ensuremath{\frac{\partial u}{\partial t}}}
\newcommand{\dbudt}{\ensuremath{\frac{\dif\bu}{\dif t}}}
\newcommand{\ptudx}{\ensuremath{\frac{\partial^2 u}{\partial x^2}}}
\newcommand{\pfudx}{\ensuremath{\frac{\partial^4 u}{\partial x^4}}}
\newcommand{\ddxi}{\ensuremath{\frac{\dif\;{}}{\dif\xi}}}
\newcommand{\dvdxi}{\ensuremath{\frac{\dif v}{\dif\xi}}}
\newcommand{\dtvdxi}{\ensuremath{\frac{\dif^2 v}{\dif\xi^2}}}
\newcommand{\dfvdxi}{\ensuremath{\frac{\dif^4 v}{\dif\xi^4}}}
\newcommand{\pdxi}{\ensuremath{\frac{\partial\;{}}{\partial\xi}}}
\newcommand{\pTdxi}{\ensuremath{\frac{\partial\cT}{\partial\xi}}}
\newcommand{\ptTdxi}{\ensuremath{\frac{\partial^2\cT}{\partial\xi^2}}}
\newcommand{\pfTdxi}{\ensuremath{\frac{\partial^4\cT}{\partial\xi^4}}}
\newcommand{\pTdp}{\ensuremath{\frac{\partial\cT}{\partial\phi}}}

\newcommand{\gam}{\ensuremath{\gamma}}
\newcommand{\lam}{\ensuremath{\lambda}}
\newcommand{\lamM}{\ensuremath{\lam_{\scriptscriptstyle M}}}
\newcommand{\lamLM}{\ensuremath{\lam_{\scriptscriptstyle L,M}}}
\newcommand{\lams}{\ensuremath{\lam^{(0)}}}
\newcommand{\lamk}{\ensuremath{\lam^{(k)}}}
\newcommand{\lamko}{\ensuremath{\lam^{(k+1)}}}
\newcommand{\lamztb}{\ensuremath{\lamz_{\mathrm{tb}}}}
\newcommand{\lamzrw}{\ensuremath{\lamz_{\mathrm{rw}}}}
\newcommand{\lamzns}{\ensuremath{\lamz_{\mathrm{ns}}}}
\newcommand{\lamz}{\ensuremath{\lam^\star}}

\newcommand{\delam}{\ensuremath{\delta\!\lam}}
\newcommand{\delamko}{\ensuremath{{\delam}^{(k+1)}}}
\newcommand{\bxzl}{\ensuremath{\bxz(\lam)}}

\newcommand{\mJxzll}{\ensuremath{\mJ(\bxz(\lam),\lam)}}
\newcommand{\mJxzlz}{\ensuremath{\mJ(\bxz,\lamz)}}
\newcommand{\mJvztlz}{\ensuremath{\mJ(\bvzt,\lamz)}}
\newcommand{\mJvztll}{\ensuremath{\mJ(\bvztl,\lam)}}
\newcommand{\mJvzttlz}{\ensuremath{\mJ(\bvzt,\theta,\lamz)}}
\newcommand{\mJvztltl}{\ensuremath{\mJ(\bvztl,\theta,\lam)}}

\newcommand{\om}{\ensuremath{\omega}}
\newcommand{\omM}{\ensuremath{\om_{\scriptscriptstyle M}}}
\newcommand{\oms}{\ensuremath{\om^{(0)}}}
\newcommand{\omk}{\ensuremath{\om^{(k)}}}

\newcommand{\omz}{\ensuremath{\om^\star}}
\newcommand{\omzns}{\ensuremath{\omz_{\mathrm{ns}}}}

\newcommand{\omtp}{\ensuremath{\om\tp}}

\newcommand{\delom}{\ensuremath{\delta\om}}
\newcommand{\delomko}{\ensuremath{{\delom}^{(k+1)}}}

\newcommand{\tpi}{\ensuremath{2\pi}}
\newcommand{\rtp}{\ensuremath{\frac{1}{\tpi}}}
\newcommand{\rtps}{\ensuremath{\frac{1}{[\tpi]^2}}}

\newcommand{\rhoLM}{\ensuremath{\rho_{\scriptscriptstyle L,M}}}
\newcommand{\rhoz}{\ensuremath{\rho^{\star}}}
\newcommand{\rhotp}{\ensuremath{\rho\tp}}

\newcommand{\Tz}{\ensuremath{T^\star}}
\newcommand{\Tzl}{\ensuremath{\Tz(\lam)}}
\newcommand{\tpTz}{\ensuremath{\tpi\Tz}}

\newcommand{\IFT}{Implicit Function Theorem}
\newcommand{\cF}{\ensuremath{\mathcal{F}}}

\newcommand{\cL}{\ensuremath{\mathcal{L}}}
\newcommand{\cT}{\ensuremath{\mathcal{T}}}
\newcommand{\Y}{\ensuremath{\mathcal{Y}}}
\newcommand{\Yp}{\ensuremath{\Y_+}}
\newcommand{\Ym}{\ensuremath{\Y_-}}
\newcommand{\Ypm}{\ensuremath{\Y_{\pm}}}
\newcommand{\Ymnm}{\ensuremath{\Ym^{\nm}}}
\newcommand{\Ypnp}{\ensuremath{\Yp^{\np}}}
\newcommand{\Ymnnm}{\ensuremath{\Ym^{n\times\nm}}}
\newcommand{\Ypnnp}{\ensuremath{\Yp^{n\times\np}}}

\title{FROM HOPF TO NEIMARK--SACKER BIFURCATION:\\ A COMPUTATIONAL ALGORITHM}
\author{Gerald Moore\thanks{Department of Mathematics, Imperial College of Science, Technology and Medicine;
	180 Queen's Gate, London SW7 2AZ (\texttt{g.moore@imperial.ac.uk})}}
\begin{document}
\maketitle
\begin{abstract}
	We construct an algorithm for approximating the invariant tori created at a Neimark--Sacker bifurcation point.
	It is based on the same philosophy as many algorithms for approximating the periodic orbits created at a Hopf bifurcation point, i.e.
	a Fourier spectral method. For Neimark--Sacker bifurcation, however, we use a simple parametrisation of the tori in order to determine low-order
	approximations, and then utilise the information contained therein to develop a more general parametrisation suitable for computing higher-order
	approximations. Different algorithms, applicable to either autonomous or periodically-forced systems of differential equations, are obtained.
\end{abstract}
\newtheorem{assumption}{Assumption}[section]
\pagestyle{myheadings}
\thispagestyle{plain}
\markboth{GERALD MOORE}{NEIMARK--SACKER BIFURCATION}
\section{Introduction}\label{Intro}
In this paper we consider both nonlinear autonomous systems
\begin{equation}\label{autde}
	\dxdt=\bFxl\qquad\bF:\Rn\times\R\mapsto\Rn,
\end{equation}
i.e. $\bF$ is a smooth function on $\Rn$ depending on a parameter $\lam$, and periodically-forced systems
\begin{equation}\label{pfde}
	\dxdt=\bFxtl\qquad\bF:\Rn\times\R\times\R\mapsto\Rn,
\end{equation}
i.e. $\bF$ also depends periodically on the independent variable $t$. In \S\ref{NSbifpf} and \S\ref{NSbifaut}, we will describe (closely related)
algorithms for approximating the invariant tori created at Neimark--Sacker bifurcation points of both \eqref{autde} and \eqref{pfde}: Neimark--Sacker
bifurcation for \eqref{pfde} is defined by Assumptions~\vrefrange{NSpffirst}{pfassum2}, while Neimark--Sacker bifurcation for \eqref{autde} is defined
by Assumptions~\vrefrange{NSautfirst}{autassum2}.
In \S\ref{Hopfbif}, however, we first introduce some of our ideas within the relatively simple paradigm case of Hopf bifurcation for \eqref{autde},
which is defined by Assumptions~\vrefrange{Hopffirst}{Hopftrans}, since the periodic orbits created here depend only on a \ul{single} frequency. In
contrast, the invariant tori created at a Neimark--Sacker bifurcation point have \ul{two} independent frequencies and it is their possible resonance that
creates the difficulties.
	
Let $(\bxz,\lamz)\in{\Rn\times\R}$ be a stationary solution of \eqref{autde}, i.e. $\bFxzlz=\zz$, at which $\mJxzlz$, the Jacobian matrix  of $\bF$,
has $n-2$ hyperbolic eigenvalues (nonzero real part) and a pair of purely imaginary eigenvalues. By the \IFT\ there is a locally unique curve
of stationary solutions, parametrised by $\lam$, satisfying
\[\bFxzll=\zz\]
and the key condition for Hopf bifurcation is Assumption~\ref{Hopftrans} on page~\pageref{Hopftrans}, that the two critical eigenvalues of $\mJxzll$
must cross the imaginary axis transversally at $\lam=\lamz$.
There is then a locally unique curve of periodic orbits for \eqref{autde} in the neighbourhood of $(\bxz,\lamz)$.
Analytical methods to investigate Hopf bifurcation are contained in \cite{All77,HKW81,HW78,Kuz04,MM76}.
In \S\ref{Hopfbif}, we show how low-order Fourier approximations of these periodic orbits simultaneously provide approximations to both the
near-identity polynomial mappings which locally transform \eqref{autde} to normal form and also to the Lyapunov coefficient in the normal form.

For Neimark--Sacker bifurcation of \eqref{autde}, we assume that $\buz(t)$ is a periodic orbit for $\lam=\lamz$, of period $\tpTz$. We also assume
that $n-3$ of the Floquet exponents of $\buz$ are hyperbolic and $2$ purely imaginary, the other being zero of course. Hence, by the \IFT, \eqref{autde}
has a locally unique curve of periodic orbits, parametrised by $\lam$, and our key condition is again Assumption~\ref{auttrans} on page~\pageref{auttrans},
i.e. that the critical pair of Floquet exponents crosses the
imaginary axis transversally at $\lam=\lamz$. In contrast to Hopf bifurcation, however, we need two additional conditions in order to guarantee the
creation of invariant tori at $(\buz,\lamz)$:
\begin{itemize}
	\item the no strong resonance Assumption~\ref{autassum1} on page~\pageref{autassum1},
		so that torus bifurcation rather than subharmonic bifurcation is generic \cite{IJ90};
	\item the real Lyapunov coefficient is nonzero in \eqref{exactLaut}, which is equivalent to the parameter $\lam$ moving away from $\lamz$ at
		leading order in Assumption~\ref{autassum2} on page~\pageref{autassum2}.
\end{itemize}
Analytical methods to investigate Neimark--Sacker bifurcation are contained in \cite{IACT81,IJ90,Wan78}. In \S\ref{secCRformaut}, we first show how
Assumption~\ref{autassum1} permits the computation of low-order Fourier approximations for our invariant tori. The information contained in these
low-order approximations is then used, together with Assumption~\ref{autassum2}, to construct higher-order Fourier approximations in \S\ref{hoFaut}.

Neimark--Sacker bifurcation of \eqref{pfde} is similar. We assume that $\buz(t)$ is a periodic orbit at $\lam=\lamz$ and also 
that $n-2$ of the Floquet exponents are hyperbolic, while $2$ are purely imaginary. Hence, by the \IFT, \eqref{pfde}
has a locally unique curve of periodic orbits, parametrised by $\lam$, and our key condition is again Assumption~\ref{pftrans} on page~\pageref{pftrans},
i.e. that the critical pair of Floquet exponents crosses the
imaginary axis transversally at $\lam=\lamz$. We still need the above two additional conditions, Assumption~\ref{pfassum1} on page~\pageref{pfassum1}
and Assumption~\ref{pfassum2} on page~\pageref{pfassum2}, in order to guarantee the creation of invariant tori at $(\buz,\lamz)$.
Analytical methods to investigate Neimark--Sacker bifurcation for \eqref{pfde} are contained in \cite{IA92,IJ90}.
In \S\ref{NSbifpf}, we again use Assumptions~\ref{pfassum1} and \ref{pfassum2} to first construct low-order and subsequently higher-order
Fourier approximations for our tori.
(We have chosen this ordering for the sections because the absence of a zero Floquet exponent simplifies our equations,
in particular the torus parametrisation is simpler. Hence transforming \eqref{pfde} to \eqref{autde}, by adding time as a new state variable, is
not recommended.)

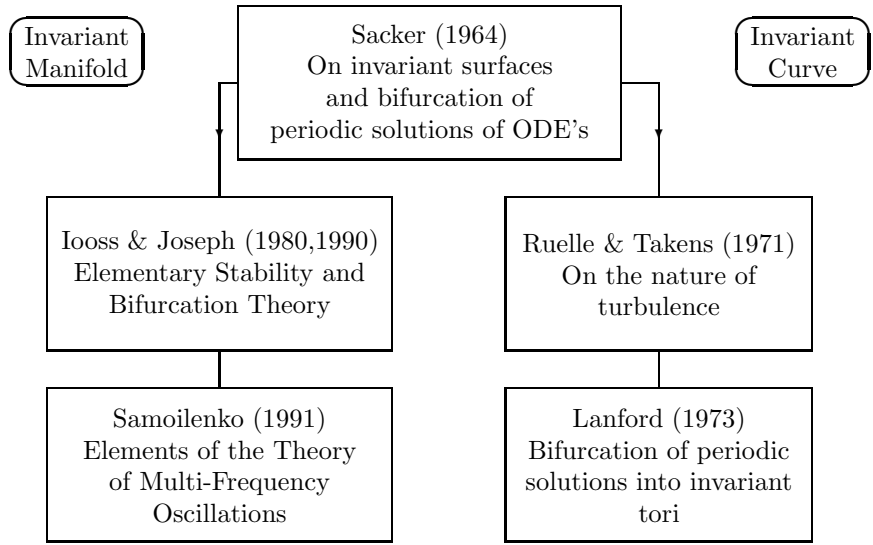
\begin{figure}[t]\setlength{\unitlength}{.04in}
	\begin{center}\begin{picture}(100,70)
		\put(0,25){\framebox(45,20){\parbox{1.8in}{\centering Iooss \& Joseph (1980,1990)\\ Elementary Stability and\\ Bifurcation Theory}}}
		\put(25,50){\framebox(50,20){\parbox{2.in}{\centering Sacker (1964)\\ On invariant surfaces\\ and bifurcation of\\
														periodic solutions of ODE's}}}
		\put(60,25){\framebox(40,20){\parbox{1.6in}{\centering Ruelle \& Takens (1971)\\ On the nature of turbulence}}}
		\put(0,0){\framebox(45,20){\parbox{1.8in}{\centering Samoilenko (1991)\\ Elements of the Theory\\ of Multi-Frequency\\ Oscillations}}}
		\put(60,0){\framebox(40,20){\parbox{1.6in}{\centering Lanford (1973)\\ Bifurcation of periodic solutions into invariant tori}}}
		\put(22.5,60){\vector(0,-1){7.5}}
		\put(22.5,52.5){\line(0,-1){7.5}}
		\put(22.5,20){\line(0,1){5}}
		\put(22.5,60){\line(1,0){2.5}}
		\put(75,60){\line(1,0){5}}
		\put(80,60){\vector(0,-1){7.5}}
		\put(80,52.5){\line(0,-1){7.5}}
		\put(80,20){\line(0,1){5}}
		\put(-5,65){\Ovalbox{\parbox[t]{.6in}{\centering Invariant\\ Manifold}}}
		\put(90,65){\Ovalbox{\parbox[t]{.6in}{\centering Invariant\\ Curve}}}
	\end{picture}\end{center}
	\caption{\label{keyref}Key references for Neimark--Sacker bifurcation}
\end{figure}

The fundamental idea behind the present paper is to use the approach in \cite{Sack64}, of which \cite{Sack65} is a concise version, to develop
the analytical foundations of a practical computational algorithm for approximating the invariant tori created at Neimark--Sacker bifurcation points.
\cite{Sack64} actually proves the existence of invariant tori in two ways:
\begin{itemize}
	\item constructing tori invariant with respect to the vector field, which is the approach used in the two key books \cite{IJ90,Sam91}
		on the left-hand side of Figure~\ref{keyref};
	\item constructing curves invariant with respect to the Poincar\'{e} map, which is the approach used in the two key papers \cite{Lan73,RT71}
		on the right-hand side of Figure~\ref{keyref}.
\end{itemize}
We do not wish to depend explicitly on the trajectories of \eqref{autde} or \eqref{pfde} and so we follow the vector field approach;
our algorithm being based on Fourier approximation and Floquet theory, in particular Floquet exponents, as introduced in \cite{Moore05}.
(Thus to appreciate fully the present paper, a familiarity with the left-hand side of Figure~\ref{keyref} is recommended.)
Hence we emphasise that in this paper our concern is with invariant tori as manifolds, and we neither consider the trajectories thereon
nor the stability of the tori. (Such questions may be answered at the post-processing stage, and are dealt with in several of the references.)
As far as we are aware, the invariant manifold approach in \cite{Sack64} has not been developed further for Neimark--Sacker bifurcation in the literature,
and has certainly not been combined with the Fourier approximation ideas in \cite{Sam91}. On the other hand, there has been quite a lot of related
work on the invariant curve approach, and we refer to \cite{Kuz04} for details and references. In the present paper, we first
see, in \S\ref{CRformHopf}, how straight-forward it is to construct Fourier approximations for the periodic orbits created at a Hopf bifurcation point,
and then attempt to generalise this algorithm for Neimark--Sacker bifurcation. The latter has two additional difficulties:
coping with possible weak resonances and implementing efficiently the \emph{ideas} behind centre manifold reduction and normal form transformation.
(Sections~\ref{Hopfbif}, \ref{NSbifpf} and \ref{NSbifaut} have been deliberately written to be as similar as possible,
both as an aid to the reader and so that the key differences stand out more clearly.)
Finally, we mention that \cite{Sam91} is not explicitly concerned with Neimark--Sacker bifurcation, merely with the continuation of invariant tori
using a Fourier-Galerkin approach. In \cite{Sack64}, however, it has already been shown how Neimark--Sacker bifurcation can be reduced to this case,
and the constructive approach in \cite{Sam91} is much more relevant to us than the uniform norm analysis based on elliptic regularisation and smoothing
operators employed in \cite{Sack64}.

\section{Hopf bifurcation}\label{Hopfbif}
We start with our two basic conditions.
\begin{assumption}\label{Hopffirst}
	$(\bxz,\lamz)\in{\Rn\times\R}$ is a stationary solution of \eqref{autde}, i.e.
	\[\bFxzlz=\zz.\]
\end{assumption}
\begin{assumption}
	$(\bxz,\lamz)$ is a Hopf bifurcation point for \eqref{autde}, i.e. $\exists\mEz\in\Rnn$ and non-singular $\mPz\in\Rnn$ such that
	\[\mJxzlz\mPz=\mPz\mEz,\]
	where
	\[\mEz\equiv\left(\begin{array}{ll}
				\mEhz&\ZZ\\
				\ZZ&\wEz
			\end{array}\right)\quad\begin{aligned}
							\mEhz&\in\Rtt\\
							\wEz&\in\Rntnt
						\end{aligned}\qquad\text{with}\quad
	\mEhz\equiv\begin{pmatrix}
			0&-\omz\\ \omz&0
			\end{pmatrix}\quad\omz>0\]
	and $\wEz$ having no eigenvalues on the imaginary axis.
\end{assumption}

\noindent This invariant subspace decomposition is a stronger assumption than required for Hopf bifurcation (the standard case of $\wEz$ having
no eigenvalue which is an integer multiple of $\ii\omz$ being considered in \cite{Moore05}) and is chosen so that this section agrees more closely with
\S\ref{NSbifpf} and \S\ref{NSbifaut}.
From these two assumptions, the \IFT\ tells us that there is a locally unique curve of stationary points, smoothly parametrised by $\lam$, and satisfying
\begin{equation}\label{Hopfcurve}
	\bFxzll=\zz.
\end{equation}
The invariant subspace decomposition may also be smoothly continued locally, and so we have
\begin{equation}\label{Hopfinvar}
	\mJxzll\mPl=\mPl\mEl\qquad{\mP:\R\mapsto\Rnn,\;\mE:\R\mapsto\Rnn},
\end{equation}
where $\mPl$ is non-singular and
\[\mEl\equiv\begin{pmatrix}
		\mEhl&\ZZ\\
		\ZZ&\wEl
	\end{pmatrix}\qquad\begin{aligned}
				\mEh&:\R\mapsto\Rtt\\
				\wE&:\R\mapsto\Rntnt
			\end{aligned}\]
with
\[\mEhl\equiv\left(\begin{array}{rr}
			\aRl&-\aIl\\ \aIl&\aRl
		\end{array}\right)\qquad\begin{aligned}
						\aR&:\R\mapsto\R\\
						\aI&:\R\mapsto\R.
					\end{aligned}\]
Finally, the key transversality condition must also hold.
\begin{assumption}\label{Hopftrans}Transversal crossing of critical eigenvalues, i.e.
	\[\boxed{\daRz\equiv\daRlz\neq0.}\]
\end{assumption}

\subsection{Crandall--Rabinowitz formulation}\label{CRformHopf}
We seek periodic orbits of \eqref{autde} in the form
\begin{equation}\label{Hopfpo}
	\bxzl+\ep\mPl\bzt\qquad\bz:\So\mapsto\Rn
\end{equation}
with unknown $\bz$, and also with unknown frequency $\om\in\R$. $\ep$ plays the role of a small amplitude parameter, upon which the unknowns $\bz$, $\lam$
and $\om$ depend. Thus the periodic orbits must satisfy
\begin{subequations}\label{HopfCR}
\begin{equation}\label{HopfCRa}
	\bF\Big(\bxzl+\ep\mPl\bzt,\lam\Big)-\om\ddt\Big\{\bxzl+\ep\mPl\bzt\Big\}=\zz
\end{equation}
and the scalar amplitude and phase conditions \cite{HKW81,IJ90,MM76,Moore05}
\begin{equation}\label{Hopfscon}
	\bgam(\bz)\equiv
		\begin{pmatrix}
			\left\langle\bzt,\bazt\right\rangle\\
			\left\langle\bzt,\dbazt\right\rangle
		\end{pmatrix}=\beh_1;
\end{equation}
\end{subequations}
here
\[\bazt\equiv(\cos{\theta},\sin{\theta},0,\dots,0)^T,\;\dbazt\equiv(-\sin{\theta},\cos{\theta},0,\dots,0)^T,\;\beh_1\equiv(1,0)^T,\]
with the inner-product defined by
\begin{equation}\label{inprod1D}
	\left\langle\bw_1(\theta),\bw_2(\theta)\right\rangle\equiv\rtp\int_0^{\tpi}\bw_1(\theta)\cdot\bw_2(\theta)\dt
\end{equation}
for $\bw_1,\bw_2:\So\mapsto\R^s$ and $s\ge1$. In order to apply the \IFT\ to \eqref{HopfCR}, we must first eliminate the curve of stationary solutions:
thus the Crandall--Rabinowitz formulation (as used in \cite{CR71} for the bifurcation of non-trivial stationary solutions) writes
\[\bG\Big(\bzt,\lam;\ep\Big)\equiv\rep\mPl^{-1}\bF\Big(\bxzl+\ep\mPl\bzt,\lam\Big)\]
and solves \eqref{HopfCR} in the form
\begin{equation}\label{CRform}
	\zz=\cF\Big(\bzt,\lam,\om;\ep\Big)\equiv
	\left\{\begin{gathered}
		\bG\Big(\bzt,\lam;\ep\Big)-\om\ddt\bzt\\
		\bgam(\bz)-\beh_1.
	\end{gathered}\right.
\end{equation}
Thus, using \eqref{Hopfcurve} and \eqref{Hopfinvar}, we can expand $\bG$ in the form
\begin{equation}\label{Gexpand}
	\bG\Big(\bz,\lam;\ep\Big)=\mEl\bz+\sum_{p\ge2}\ep^{p-1}\bG_p\Big(\bz;\lam\Big)\qquad\bG_p:\Rn\times\R\mapsto\Rn,
\end{equation}
the $n$ components of $\bG_p$ being homogeneous polynomials of degree $p$ in the $n$ components of $\bz$ with coefficients depending on $\lam$.
(Here and later we display important functions and mappings in this way; with the understanding that the sum is limited by the available smoothness.) 
At $\ep=0$, \eqref{CRform} has the solution 
\[\bzt=\bazt,\;\lam=\lamz,\;\om=\omz\]
and the linearisation of \eqref{CRform} about this solution is \ul{non-singular} since a simple Fourier analysis (using the properties of $\mEz$)
shows that
\begin{gather*}
	\left[\mEz-\omz\mI\ddt\right]\bzt+\lam\dElz\bazt-\om\dbazt=\bft\\
	\bgam(\bz)=\zz
\end{gather*}
implies the existence of a constant $C_{\cL}>0$ such that
\[\max\left\{\lVert\bz\rVert_{H^1},\lvert\lam\rvert,\lvert\om\rvert\right\}\le C_{\cL}\lVert\bff\rVert_{L^2}.\]
(Here we use the standard spaces/norms of periodic functions \cite{Sam91}, based on the inner-product \eqref{inprod1D}.)
Hence the \IFT, relying on a Newton-chord iteration for constructing solutions of \eqref{CRform} from the starting value
\[\bzst=\bazt,\;\lams=\lamz,\;\oms=\omz,\]
gives the following result.
\begin{theorem}\label{HopfEU}
	For all $\lvert\ep\rvert$ sufficiently small, \eqref{CRform} has a locally unique solution
	\[\lamz(\ep),\quad\omz(\ep)\quad\text{and}\quad\bzz\tep.\]
	It can be written as an expansion in powers of $\ep$ \cite{IJ90}, i.e. 
	\begin{equation}\label{Hopfsol}
		\begin{gathered}
			\lamz(\ep)\equiv\lamz+\sum_{p\ge1}\ep^{2p}\lamz_{2p},\quad\omz(\ep)\equiv\omz+\sum_{p\ge1}\ep^{2p}\omz_{2p},\\
			\bzz\tep\equiv\bazt+\sum_{p\ge1}\ep^{2p-1}\bzz_{2p}(\theta)+\ep^{2p}\bzz_{2p+1}(\theta);
		\end{gathered}
	\end{equation}
	where $\bzz_{2p}$ only depends on the \ul{even} Fourier modes $0,2,4,\dots,2p$ and $\bzz_{2p+1}$ only depends on the \ul{odd} Fourier modes
	$1,3,5,\dots,2p+1$. The amplitude and phase conditions force
	\[\left\langle\bzz_{2p+1}(\theta),\bazt\right\rangle=0=\left\langle\bzz_{2p+1}(\theta),\dbazt\right\rangle.\]
	\eqref{Hopfsol} can also be expressed in Fourier modes, i.e.
	\begin{equation}\label{HopfsolF}
		\bzz\tep\equiv\bazt+\baz_0(\ep)+\sum_{m\ge1}\baz_m(\ep)\cos{m\theta}+\bbz_m(\ep)\sin{m\theta},
	\end{equation}
	where
	\[\begin{array}{crcl}
		&\baz_0(\ep)&\text{has terms}&\ep,\ep^3,\ep^5,\dots\\
		&\baz_1(\ep),\bbz_1(\ep)&\text{have terms}&\ep^2,\ep^4,\ep^6,\dots\\
		m\ge2&\baz_m(\ep),\bbz_m(\ep)&\text{have terms}&\ep^{m-1},\ep^{m+1},\ep^{m+3},\dots.
	\end{array}\]
	Again, the amplitude and phase conditions force
	\begin{equation}\label{Hopfap}
		\baz_1(\ep)\cdot\bef+\bbz_1(\ep)\cdot\bet=0\quad\text{and}\quad\baz_1(\ep)\cdot\bet-\bbz_1(\ep)\cdot\bef=0.
	\end{equation}
\end{theorem}

In practice we can construct accurate approximations to our periodic orbits by computing $\lam$, $\om$ and a finite Fourier series
\[\bzMt\equiv\bazt+\ba_0+\sum_{m=1}^M\ba_m\cos{m\theta}+\bb_m\sin{m\theta}\]
which solve the Galerkin equations for \eqref{CRform}; i.e.
\begin{equation}\label{CRM}
\begin{gathered}
	\mSM\bG\Big(\bzMt,\lam;\ep\Big)-\om\ddt\bzMt=\zz\\
	\bgam(\bzM)-\beh_1=\zz,
\end{gathered}
\end{equation}
where $\mSM:{L^2\mapsto L^2}$ is the operator which performs the Fourier series truncation.
Thus we have the usual approximation result in terms of the decay of the Fourier modes in \eqref{HopfsolF}.
\begin{theorem}\label{HopfFapprox}
	For all $\lvert\ep\rvert$ sufficiently small, \eqref{CRM} has a locally unique solution
	\begin{equation}\label{HopfM}
		\begin{gathered}
			\lamM^F(\ep),\;\omM^F(\ep)\;\text{and}\\
			\bzM^F\tep\equiv\bazt+\ba_0^F(\ep)+\sum_{m=1}^M\ba_m^F(\ep)\cos{m\theta}+\bb_m^F(\ep)\sin{m\theta},
		\end{gathered}
	\end{equation}
	which satisfies the error bounds
	\begin{multline*}
		\max\left\{\lVert\bzM^F(.;\ep)-\mSM\bzz(.;\ep)\rVert_{H^1},
			\lvert\lamM^F(\ep)-\lamz(\ep)\rvert,\lvert\omM^F(\ep)-\omz(\ep)\rvert\right\}\\
		\le C_{\cF}\lVert\left(\mI-\mSM\right)\bzz(.;\ep)\rVert_{H^1}.
	\end{multline*}
\end{theorem}
(In this paper, we will not be considering any superconvergence phenomena.)

The Fourier approximation in Theorem~\ref{HopfFapprox} has no restriction on the size of $M$ and, similarly to \eqref{Hopfsol},
it can also be written as an asymptotic expansion in powers of $\ep$. Thus instead of considering the approximation error for fixed $\ep$ as $M$ increases,
it is also possible to consider this error for fixed small $M$ as $\ep\to0$. In \S\ref{approxnf}, we will
make particular use of the approximation for \ul{$M=3$}, i.e.
\begin{equation}\label{HopfF3}
	\begin{gathered}
		\lam_3^F(\ep),\;\om_3^F(\ep)\;\text{and}\\
		\bz_3^F\tep\equiv\bazt+\ba_0^F(\ep)+\sum_{m=1}^3\ba_m^F(\ep)\cos{m\theta}+\bb_m^F(\ep)\sin{m\theta},
	\end{gathered}
\end{equation}
where, as in \eqref{Hopfap}, the amplitude and phase conditions force
\begin{equation}\label{Fourap}
	\ba_1^F(\ep)\cdot\bef+\bb_1^F(\ep)\cdot\bet=0\quad\text{and}\quad\ba_1^F(\ep)\cdot\bet-\bb_1^F(\ep)\cdot\bef=0.
\end{equation}
\begin{corollary}\label{HopfF3err}
	Comparing \eqref{HopfF3} with the exact solution in \eqref{HopfsolF} gives the errors
	\begin{align*}
		&&\left\lvert\lamz(\ep)-\lam_3^F(\ep)\right\rvert&=O(\ep^4),&\left\lvert\omz(\ep)-\om_3^F(\ep)\right\rvert&=O(\ep^4),&\\
		m&=0,2&\left\lVert\baz_m(\ep)-\ba_m^F(\ep)\right\rVert&=O(\ep^3),&\left\lVert\bbz_m(\ep)-\bb_m^F(\ep)\right\rVert&=O(\ep^3)&m&=2,\\
		m&=1,3&\left\lVert\baz_m(\ep)-\ba_m^F(\ep)\right\rVert&=O(\ep^4),&\left\lVert\bbz_m(\ep)-\bb_m^F(\ep)\right\rVert&=O(\ep^4)&m&=1,3.
	\end{align*}
\end{corollary}

\subsection{Normal form and its Fourier approximation}\label{approxnf}
Our algorithm for Hopf bifurcation in \S\ref{CRformHopf} requires neither reduction to the centre manifold nor transformation to normal form.
For Neimark--Sacker bifurcation, however, these two procedures have to be implemented approximately in order to cope properly with possible weak resonances.
Thus we now choose to illustrate our later approach in the present relatively simple setting.

Instead of carrying out the standard theoretical centre manifold reduction and normal form tranformation \cite{GH83,Kuz04,Mur03}, we adopt the
\ul{operational} approach in \cite{CS83,IA92,IJ90,Kuz04} and construct the necessary transformations in order to simplify the key equation \eqref{CRform},
i.e.
\begin{equation}\label{cmnfde}
	\begin{gathered}
		\bG\Big(\bzt,\lam;\ep\Big)-\om\ddt\bzt=\zz\\
		\bgam(\bz)-\beh_1=\zz.
	\end{gathered}
\end{equation}
By introducing
\[\bz^T\equiv\left(\bzh^T,\wtz^T\right),\quad\text{with $\bzh\in\Rt$ and $\wtz\in\Rnt$,}\]
we first write
\[\bG\Big(\bz,\lam;\ep\Big)\equiv
	\begin{bmatrix}
		\bGh\Big(\bzh,\wtz,\lam;\ep\Big)\\
		\wbG\Big(\bzh,\wtz,\lam;\ep\Big)
	\end{bmatrix}\qquad
	\begin{aligned}
		\bGh:{\Rt\times\Rnt\times\R\times\R}&\mapsto\Rt\\
		\wbG:{\Rt\times\Rnt\times\R\times\R}&\mapsto\Rnt.
	\end{aligned}\]
We then aim to simplify the lower terms in $\bGh$ and $\wbG$ as much as possible by constructing suitable mappings
\[\bhh:{\Rt\times\R}\mapsto\Rt\quad\text{and}\quad\wbh:{\Rt\times\R}\mapsto\Rnt,\]
where $\wbh$ is a homogeneous quadratic polynomial with $\lam$-dependent coefficients and
$\bhh$ is the sum of homogeneous quadratic and cubic polynomials with $\lam$-dependent coefficients, to define the near-identity transformations
\begin{subequations}\label{hhtdef}
\begin{equation}
\begin{aligned}
	\bzh&=\byh+\rep\bhh(\ep\byh;\lam)\\
	&=\byh+\ep\left\{\left[\hy_1^2+\hy_2^2\right]\bah_0(\lam)+\left[\hy_1^2-\hy_2^2\right]\bah_2(\lam)+2\hy_1\hy_2\bbh_2(\lam)\right\}\\
	&\qquad+\ep^2\left\{\hy_1\left[\hy_1^2+\hy_2^2\right]\bah_1(\lam)+\hy_2\left[\hy_1^2+\hy_2^2\right]\bbh_1(\lam)\right.\\
		&\qquad\qquad\left.+\hy_1\left[\hy_1^2-3\hy_2^2\right]\bah_3(\lam)+\hy_2\left[3\hy_1^2-\hy_2^2\right]\bbh_3(\lam)\right\}
\end{aligned}
\end{equation}
and
\begin{equation}
\begin{aligned}
	\wtz&=\wty+\rep\wbh(\ep\byh;\lam)\\
	&=\wty+\ep\left\{\left[\hy_1^2+\hy_2^2\right]\wba_0(\lam)+\left[\hy_1^2-\hy_2^2\right]\wba_2(\lam)+2\hy_1\hy_2\wbb_2(\lam)\right\}.
\end{aligned}
\end{equation}
\end{subequations}
The homogeneous polynomials are given the above bases in order to link up with the Fourier coefficients through elementary trigonometrical identities,
as the table in Figure~\ref{modepoly} shows.
\begin{figure}[h!]\begin{center}
	\begin{minipage}{\textwidth}
		\[\begin{array}{|c|c||c|c||c|c|}\hline
			\cos\theta	&	y_1&
			\cos{2\theta}	&	y_1^2-y_2^2&
			\cos{3\theta}	&	y_1(y_1^2-3y_2^2)\\ \hline
			\sin\theta	&	y_2&
			\sin{2\theta}	&	2y_1y_2&
			\sin{3\theta}	&	y_2(3y_1^2-y_2^2)\\ \hline
		\end{array}\]
	\end{minipage}\end{center}
	\caption{\label{modepoly}Linking Fourier modes and polynomials}
\end{figure}
(Of course, by writing our Fourier series in exponential form, this correspondence is simpler; but we do not wish to give the impression that complex
arithmetic is necessary.) Thus we see how (through $y_1^2+y_2^2=1$)
the \ul{resonant} cubic terms, the null-space of the adjoint of the homological operator in the usual normal form computations
\cite{IA92,Mur03} being spanned by
\begin{equation}\label{cubres}
	\left[y_1^2+y_2^2\right]\begin{pmatrix}y_1\\y_2\end{pmatrix}\quad\text{and}\quad
		\left[y_1^2+y_2^2\right]\left(\begin{array}{r}-y_2\\y_1\end{array}\right),
\end{equation}
appear through these identities, and how we must have the restrictions 
\begin{equation}
	\bah_1(\lam)\cdot\beh_1+\bbh_1(\lam)\cdot\beh_2=0\quad\text{and}\quad\bah_1(\lam)\cdot\beh_2-\bbh_1(\lam)\cdot\beh_1=0
\end{equation}
in the definition of $\bhh$. Under these near-identity tranformations, \eqref{cmnfde} becomes
\begin{subequations}\label{csubde}
\begin{gather}
	\bGhd\Big(\byht,\wtyt,\lam;\ep\Big)-\om\ddt\byht=\zz\label{csubdea}\\
	\wbGd\Big(\byht,\wtyt,\lam;\ep\Big)-\om\ddt\wtyt=\zz\label{csubdeb}\\
	\bgam(\byh)-\beh_1=\zz:\label{csubdec}
\end{gather}
\end{subequations}
the two mappings
\[\bGhd:{\Rt\times\Rnt\times\R\times\R}\mapsto\Rt\quad\text{and}\quad\wbGd:{\Rt\times\Rnt\times\R\times\R}\mapsto\Rnt\]
capable of being expanded, like \eqref{Gexpand}, in the form
\begin{align*}
	\bGhd\Big(\byh,\wty,\lam;\ep\Big)&=\mEhl\byh+\sum_{p\ge2}\ep^{p-1}\bGhd_p\Big(\byh,\wty;\lam\Big)\qquad\bGhd_p:\Rt\times\Rnt\times\R\mapsto\Rt\\
	\wbGd\Big(\byh,\wty,\lam;\ep\Big)&=\wEl\wty+\sum_{p\ge2}\ep^{p-1}\wbGd_p\Big(\byh,\wty;\lam\Big)\qquad\wbGd_p:\Rt\times\Rnt\times\R\mapsto\Rnt,
\end{align*}
where the components of $\bGhd_p$ and $\wbGd_p$ are homogeneous polynomials of degree $p$ in the components of $\byh$ and $\wty$, with coefficients
depending on $\lam$. Now we choose $\bhh$ and $\wbh$ so that the lower terms in $\bGhd$ and $\wbGd$ may be simplified in the following way:
\begin{itemize}
	\item $\wbh$ forces the coefficients of the quadratic terms for $\byh$ in $\wbGd_2$ to be zero
	\item $\bhh$ forces the coefficients of the quadratic terms for $\byh$ in $\bGhd_2$ to be zero
		and the coefficients of the cubic terms for $\byh$ in $\bGhd_3$ to take the form
		\begin{equation}\label{LyapHopf}
			\left[\hy_1^2+\hy_2^2\right]\mBhl\byh,\qquad\text{where}\quad
			\mBhl\equiv\left(\begin{array}{rr} \bRl&-\bIl\\ \bIl&\bRl \end{array}\right)
		\end{equation}
		and we call the elements of this matrix Lyapunov coefficients.
\end{itemize}
(I.e. after transformation, only a multiple of the resonant cubic terms \eqref{cubres} remains.)

After this simplification, if we now insert
\begin{subequations}\label{csubsol}
\begin{gather}
	\byht=\bahzt\;\left(\equiv(\cos{\theta},\,\sin{\theta})^T\right)\quad\text{and}\quad\wtyt=\zz\label{csubsola}\\
	\lam=\lamz-\ep^2\frac{\bRlz}{\daRlz}\quad\text{and}\quad\om=\omz+\ep^2\frac{\daRlz\bIlz-\daIlz\bRlz}{\daRlz}\label{csubsolb}
\end{gather}
\end{subequations}
into the left-hand side of \eqref{csubde}, we can easily see that the remainder is $O(\ep^3)$ for \eqref{csubdea}, $O(\ep^2)$ for \eqref{csubdeb} and
zero for \eqref{csubdec}. Consequently, by transforming \eqref{csubsola} back through \eqref{hhtdef}, i.e.
\begin{equation}\label{csubsoltrans}
	\bzht=\bahzt+\rep\bhh(\ep\bahzt;\lam)\quad\text{and}\quad\wtzt=\rep\wbh(\ep\bahzt;\lam),
\end{equation}
we obtain an asymptotic solution for \eqref{cmnfde}. Since Theorem~\ref{HopfEU} already displays such a solution, i.e. $\lamz(\ep)$, $\omz(\ep)$ and
\[\bzz\tep\equiv\begin{pmatrix}\bzhz\tep\\ \wtzz\tep\end{pmatrix},\]
this must match with \eqref{csubsolb} and \eqref{csubsoltrans}.
Thus we obtain
\begin{equation}\label{cmnfsola}
\begin{aligned}
	\lamz(\ep)&=\lamz-\ep^2\frac{\bRlz}{\daRlz}+O(\ep^4)\\
	\omz(\ep)&=\omz+\ep^2\frac{\daRlz\bIlz-\daIlz\bRlz}{\daRlz}+O(\ep^4)
\end{aligned}
\end{equation}
and, through \eqref{csubsoltrans},
\begin{subequations}\label{cmnfsolc}
\begin{align}
	\bzhz\tep&=\bahzt+\ep\left\{\bah_0(\lamz)+\bah_2(\lamz)\cos{2\theta}+\bbh_2(\lamz)\sin{2\theta}\right\}\notag\\
		&\qquad+\ep^2\left\{\bah_1(\lamz)\cos{\theta}+\bbh_1(\lamz)\sin{\theta}\right.\notag\\
		&\qquad\qquad\left.+\bah_3(\lamz)\cos{3\theta}+\bbh_3(\lamz)\sin{3\theta}\right\}+O(\ep^3)\\
	\wtzz\tep&=\ep\left\{\wba_0(\lamz)+\wba_2(\lamz)\cos{2\theta}+\wbb_2(\lamz)\sin{2\theta}\right\}+O(\ep^2).
\end{align}
\end{subequations}
Finally, by comparing \eqref{cmnfsolc} and \eqref{Hopfsol}, we see that the coefficients defining $\bhh(.;\lamz)$ and $\wbh(.;\lamz)$ in \eqref{hhtdef}
are given \ul{exactly} by the coefficients of the Fourier modes in the
$\bzhz_2(\theta)$ and $\bzhz_3(\theta)$ terms of $\bzhz\tep$ and in the $\wtzz_2(\theta)$ term of $\wtzz\tep$ for \eqref{Hopfsol}.
Moreover, by comparing \eqref{cmnfsola} and \eqref{Hopfsol}, we also see that the Lyapunov coefficients $\bRlz$ and $\bIlz$ in \eqref{LyapHopf}
are given \ul{exactly} by
\[\bRlz=-\daRlz\lamz_2\quad\text{and}\quad\bIlz=\omz_2-\daIlz\lamz_2.\]

To calculate the expansion in \eqref{Hopfsol}, however, requires (through $\bG$) explicit knowledge of the second and third derivatives of $\bF$,
so it is practically much more convenient to approximate not only the coefficients of $\bhh(.;\lamz)$ and $\wbh(.;\lamz)$ but also the Lyapunov coefficients
$\bRlz$ and $\bIlz$ by using instead the \ul{$M=3$} Fourier approximation in \eqref{HopfF3}, i.e.
$\lam_3^F(\ep)$, $\om_3^F(\ep)$ and
\[\bz_3^F\tep\equiv\left\{
\begin{aligned}
	\bzh_3^F\tep&\equiv\bah_0^F(\ep)+\sum_{m=1}^3\bah_m^F(\ep)\cos{m\theta}+\bbh_m^F(\ep)\sin{m\theta}\\
	\wtz_3^F\tep&\equiv\wba_0^F(\ep)+\sum_{m=1}^3\wba_m^F(\ep)\cos{m\theta}+\wbb_m^F(\ep)\sin{m\theta}.
\end{aligned}
\right.\]
\begin{theorem}\label{Hopfpracerr}
	Using the asymptotic error results in Corollary~\ref{HopfF3err} on page~\pageref{HopfF3err}, our practical approximate formulae become
	\begin{align*}
		\bRlz&=-\daRlz\frac{\lam_3^F(\ep)-\lamz}{\ep^2}+O(\ep^2),\\
		\bIlz&=\frac{\om_3^F(\ep)-\omz}{\ep^2}-\daIlz\frac{\lam_3^F(\ep)-\lamz}{\ep^2}+O(\ep^2),
	\end{align*} 
	and
	{\renewcommand{\arraystretch}{1.5}
	\[\begin{array}{|l||l|l|}\hline
		m=0,2&\bah_m(\lamz)=\rep\bah_m^F(\ep)+O(\ep^2)&\wba_m(\lamz)=\rep\wba_m^F(\ep)+O(\ep^2)\\ \hline
		m=0&\bbh_m(\lamz)=\rep\bbh_m^F(\ep)+O(\ep^2)&\wbb_m(\lamz)=\rep\wbb_m^F(\ep)+O(\ep^2)\\ \hline
		m=1,3&\bah_m(\lamz)=\rept\bah_m^F(\ep)+O(\ep^2)&\bbh_m(\lamz)=\rept\bbh_m^F(\ep)+O(\ep^2)\\ \hline
	\end{array}\]}
\end{theorem}

We conclude by emphasizing how the $M=3$ Fourier results will be used later in Neimark--Sacker bifurcation. For a chosen value of $\ep$,
we can easily compute $\bz_3^F\tep$, $\lam_3^F(\ep)$ and $\om_3^F(\ep)$ from Theorem~\ref{HopfFapprox}:
the two scalar outputs then give us approximations for the Lyapunov coefficients $\bRlz$ and $\bIlz$,
while the Fourier components of $\bz_3^F\tep$ provide approximations for the coefficients of the polynomials $\bhh(.;\lamz)$ and $\wbh(.;\lamz)$.
With regard to Hopf bifurcation itself, the above approximate formulae may be regarded as alternatives to those suggested in \cite{HKW81,Kuz04,MM76}.

\subsection{Numerical results}
\begin{figure}[t]
	\begin{center}
		\includegraphics[width=.9\textwidth]{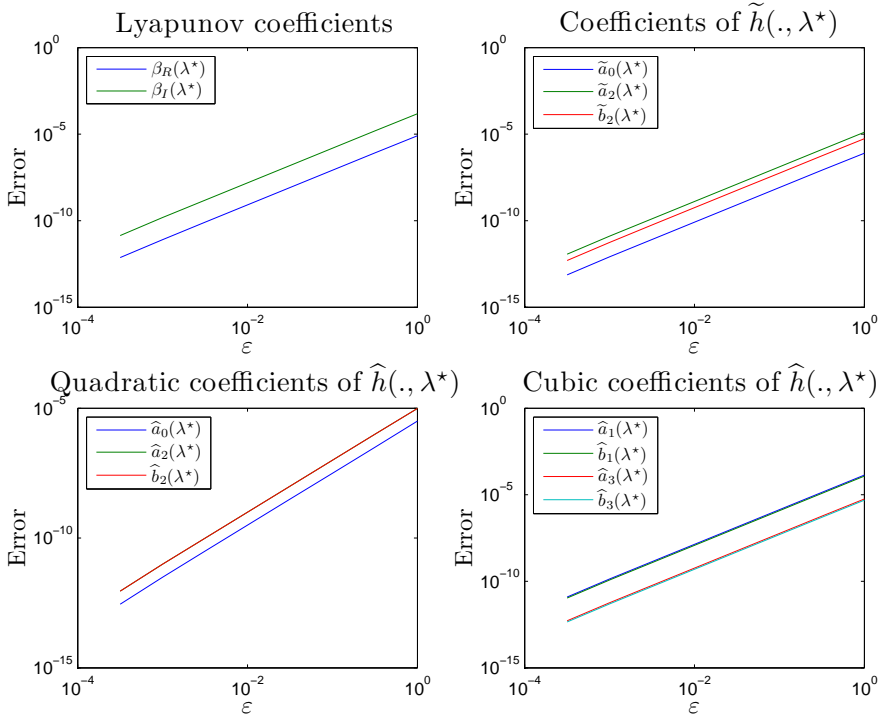}
	\end{center}
	\caption{Errors in approximation from Theorem~\ref{Hopfpracerr}\label{Hopffig}}
\end{figure}

Now we illustrate the above approximations on the well-known Lorenz equations \cite{GH83,HKW81,Kuz04,MM76}
\[\dot{x}_1=\sigma(x_2-x_1)\quad\dot{x}_2=\lam x_1-x_2-x_1x_3\quad\dot{x}_3=x_1x_2-bx_3.\]
$\sigma$ and $b$ are regarded as fixed parameters and $\lam$ is our continuation parameter.
For $\sigma>b+1$ there is a subcritical Hopf bifurcation from the stationary solution curve
\begin{gather*}
\text{$\bxzl\equiv\left(\sqrt{b(\lam-1)},\sqrt{b(\lam-1)},\lam-1\right)^T$ for $\lam>1$ at $\lamz\equiv\frac{\sigma(\sigma+b+3)}{\sigma-b-1}$,}\\
	\text{where $\omz\equiv\sqrt{b(\lamz+\sigma)}$, $\wEz\equiv-(\sigma+b+1)$ and $\daRz\equiv \frac{b(\sigma-b-1)}{2({\omz}^2+(\sigma+b+1)^2)}$.}
\end{gather*}
We use the standard parameter values $(\sigma,b)=(10,\frac83)$, which gives $\lamz\approx 24.74$, and Figure~\ref{Hopffig} displays the error for
the approximations contained in Theorem~\ref{Hopfpracerr}. Thus the $O(\ep^2)$ convergence is verified.

\section{Computational Floquet Theory}\label{CompFT}
Floquet theory enables us to transform linear, periodic ode's to constant-coefficient form: this both simplifies the analysis and leads to much more
efficient approximation by Fourier methods. A detailed discussion is contained in \cite{Moore05},
here we only describe concisely the results that are required. If the linear, periodic system we wish to solve is
\begin{equation}\label{lpode}
	-\dvdt(\theta)+\mAt\bvt=\bft\qquad\bv,\bff:\So\mapsto\Rn,\;\mA:\So\mapsto\Rnn,
\end{equation}
then our Floquet-values and Floquet-vectors solve the corresponding eigen-problem
\begin{equation}\label{Floquev}
	-\dPt+\mAt\mPt=\mPt\mE\qquad\mE\in\Rnn,\;\mPt\in\Rnn.
\end{equation}
In general, to avoid the explicit use of complex arithmetic, it is necessary to work with both periodic $(\Yp)$ and anti-periodic $(\Ym)$ mappings, i.e.
if $s:{\R\mapsto\R}$ then
\[s\in\Ypm\quad\text{iff}\quad s(\theta+\tpi)=\pm s(\theta)\;\forall\theta\in\R.\]
Thus, more specifically, $\mE$ and $\mP$ in \eqref{Floquev} have the form 
\begin{align*}
	\mE&\equiv\begin{pmatrix} \mEp&\ZZ\\ \ZZ&\mEm\end{pmatrix}&\mEp&\in\Rnpnp,\;\mEm\in\Rnmnm\\
	\mP&\equiv\begin{pmatrix} \mPp&\mPm\end{pmatrix}&\mPp&\in\Ypnnp,\;\mPm\in\Ymnnm,
\end{align*}
i.e.
\[\mPpmt\in\Rnnpm\quad\text{and}\quad\mPpmtpt=\pm\mPpmt\qquad\forall\theta\in\R,\]
for some $\np,\nm\in\Z$ with $\np+\nm=n$ and $0\leq\np\leq n$.
Then to transform \eqref{lpode} to constant-coefficient form, we transform $\bv,\bff$ to Floquet variables
\begin{align*}
	\bvt=\mPt\bwt&\equiv\begin{bmatrix}\mPpt\bwpt\\ \mPmt\bwmt\end{bmatrix}&\bwp&\in\Ypnp,\;\bwm\in\Ymnm\\
	\bft=\mPt\bgt&\equiv\begin{bmatrix}\mPpt\bgpt\\ \mPmt\bgmt\end{bmatrix}&\bgp&\in\Ypnp,\;\bgm\in\Ymnm
\end{align*}
and hence arrive at the two equations
\[-\ddt\bwpt+\mEp\bwpt=\bgpt\qquad -\ddt\bwmt+\mEm\bwmt=\bgmt\]
in $\Ypnp$ and $\Ymnm$ respectively.

Finally, we emphasise that $\npm$ are not in general unique, but can always be chosen so that the imaginary parts of the eigenvalues of $\mEpm$
(\emph{the Floquet exponents}) lie in $(-\thalf,\thalf)$.
\[\fbox{\parbox{0.66\textwidth}{
	\textbf{N.B.} \ul{For simplicity}, we shall assume in \S\ref{NSbifpf} and \S\ref{NSbifaut} that $\np=n$ and $\nm=0$.
	This is briefly commented on in \S\ref{Conclu}.}}\]

\section{Neimark--Sacker bifurcation for periodically-forced systems}\label{NSbifpf}
We may assume that the forcing in \eqref{pfde} is $\tpi$-periodic, and emphasise this by using $\theta$ as the independent variable from now on, i.e.
\eqref{pfde} becomes
\begin{equation}\label{pfdea}
	\dvdt(\theta)=\bFvtl\qquad\bF:\Rn\times\So\times\R\mapsto\Rn.
\end{equation}
We start with our two basic conditions.
\begin{assumption}\label{NSpffirst}
	At $\lam=\lamz$, $\bvz:\So\mapsto\Rn$ is a periodic orbit of \eqref{pfdea}, i.e.
	\[\bFvzttlz-\ddt\bvzt=\zz\qquad\forall\theta\in\So.\]
\end{assumption}
\begin{assumption}
	If we apply the Floquet theory in \S\ref{CompFT} to
	\[\mAt\equiv\mJvzttlz\]
	then \eqref{Floquev} becomes
	\[\mJvzttlz\mPzt-\dPzt=\mPzt\mEz,\]
	where $\mEz\in\Rnn$ and $\mPz:\So\mapsto\Rnn$ with $\mPzt$ non-singular $\forall\theta\in\So$, and we have the invariant subspace decomposition
	\[\mEz\equiv\left(\begin{array}{ll}
				\mEhz&\ZZ\\
				\ZZ&\wEz
			\end{array}\right)\quad\begin{aligned}
							\mEhz&\in\Rtt\\
							\wEz&\in\Rntnt
						\end{aligned}
	\qquad\text{with}\quad
	\mEhz\equiv\begin{pmatrix}
			0&-\omz\\ \omz&0
			\end{pmatrix}\quad\omz>0\]
	and $\wEz$ having no eigenvalues on the imaginary axis.
\end{assumption}

\noindent The \IFT\ then gives us a locally unique curve of periodic orbits, smoothly parametrised by $\lam$,
and satisfying
\begin{equation}\label{pfcurve}
	\bFvztltl-\ddt\bvztl=\zz.
\end{equation}
The Floquet variables in the invariant subspace decomposition can also be smoothly continued locally, and so we have
\begin{equation}\label{pfinvar}
	\mJvztltl\mPtl-\dPtl=\mPtl\mEl\qquad\begin{gathered}\mP:{{\So\times\R}\mapsto\Rnn}\\\mE:{\R\mapsto\Rnn},\end{gathered}
\end{equation}
where $\mPtl$ is non-singular and
\[\mEl\equiv\begin{pmatrix}
		\mEhl&\ZZ\\
		\ZZ&\wEl
	\end{pmatrix}\qquad\begin{aligned}
				\mEh&:\R\mapsto\Rtt\\
				\wE&:\R\mapsto\Rntnt
			\end{aligned}\]
with
\[\mEhl\equiv\left(\begin{array}{rr}
			\aRl&-\aIl\\ \aIl&\aRl
		\end{array}\right)\qquad\begin{aligned}
						\aR&:\R\mapsto\R\\
						\aI&:\R\mapsto\R.
					\end{aligned}\]
Finally, the key transversality condition must also hold.
\begin{assumption}\label{pftrans}Transversal crossing of critical Floquet exponents, i.e.
	\[\boxed{\daRz\equiv\daRlz\neq0.}\]
\end{assumption}
	
\subsection{Crandall--Rabinowitz formulation}\label{secCRformpf}
To start with, we attempt to mimic our approach for Hopf bifurcation in \S\ref{CRformHopf} and seek invariant tori of \eqref{pfdea} in the form
\begin{equation}\label{toruspf1}
	\bvztl+\ep\mPtl\bztp\qquad\bz:\So\times\So\mapsto\Rn,
\end{equation}
with unknown $\bz$, satisfying
\begin{subequations}\label{invarpf1}
\begin{multline}\label{invarpf1a}
	\bF\Big(\bvztl+\ep\mPtl\bztp,\theta,\lam\Big)\\
	-\pdt\Big[\bvztl+\ep\mPtl\bztp\Big]\\
	-\om\pdp\Big[\bvztl+\ep\mPtl\bztp\Big]=\zz
\end{multline}
for some unknown $\om\in\R$.
Thus we are no longer following trajectories of \eqref{pfde}, but characterising the invariance of \eqref{toruspf1} by insisting that the
vector field $\bf$ must lie in its tangent space \cite{Moore95b,Moore96}. \eqref{toruspf1} and \eqref{invarpf1a} are based on a particularly simple
choice of parametrisation for our torus, and we shall see in \S\ref{hoFpf} that more subtlety is required later. The present choice, however, is the natural
analogue of Hopf bifurcation (with $\om$ playing the role of frequency in $\phi$) and enables us to approximate the normal form in \S\ref{approxnfpf}.
Of course, we also require the scalar amplitude and phase conditions
\begin{equation}\label{pfscon}
	\bgam(\bz)\equiv
		\begin{pmatrix}
			\left\langle\left\langle\bztp,\bazp\right\rangle\right\rangle\\
			\left\langle\left\langle\bztp,\dbazp\right\rangle\right\rangle
		\end{pmatrix}=\beh_1,
\end{equation}
\end{subequations}
where
\[\bazp\equiv(\cos{\phi},\sin{\phi},0,\dots,0)^T,\quad\dbazp\equiv(-\sin{\phi},\cos{\phi},0,\dots,0)^T\]
with the inner-product defined by
\[\left\langle\left\langle\bw_1\tp,\bw_2\tp\right\rangle\right\rangle
	\equiv\rtps\int_0^{\tpi}\int_0^{\tpi}\bw_1\tp\cdot\bw_2\tp\dt\dphi.\]

To attempt to apply the \IFT\ to \eqref{invarpf1}, we must first eliminate the curve of periodic orbits: thus the Crandall--Rabinowitz formulation writes
\begin{multline*}
	\bG\Big(\bztp,\lam,\theta;\ep\Big)\equiv\rep\mPtl^{-1}\bigg\{\bF\Big(\bvztl+\ep\mPtl\bztp,\theta,\lam\Big)\\
	-\bF\Big(\bvztl,\theta,\lam\Big)-\ep\dPtl\bztp\bigg\}
\end{multline*}
and solves \eqref{invarpf1} in the form
\begin{equation}\label{CRformpf}
	\zz=\cF\Big(\bztp,\lam,\om,\theta;\ep\Big)\equiv
	\left\{\begin{gathered}
		\bG\Big(\bztp,\lam,\theta;\ep\Big)-\pbzdttp-\om\pbzdptp\\
		\bgam(\bz)-\beh_1.
	\end{gathered}\right.
\end{equation}
Hence, using \eqref{pfcurve} and \eqref{pfinvar}, we can expand $\bG$ in the form
\begin{equation}\label{Gexpandpf}
	\bG\Big(\bz,\lam,\theta;\ep\Big)=\mEl\bz+\sum_{p\ge2}\ep^{p-1}\bG_p\Big(\bz;\lam,\theta\Big)\qquad\bG_p:\Rn\times\R\times\So\mapsto\Rn,
\end{equation}
the $n$ components of $\bG_p$ being homogeneous polynomials of degree $p$ in the $n$ components of $\bz$ with coefficients depending on $\lam$ and $\theta$.
At $\ep=0$, \eqref{CRformpf} has the solution 
\[\bztp=\bazp,\;\lam=\lamz,\;\om=\omz\]
and the linearisation about this solution is
\begin{equation}\label{toruspflin}
\begin{gathered}
	\left[\mEz-\mI\pdt-\omz\mI\pdp\right]\bztp+\lam\dElz\bazp-\om\dbazp\\
	\bgam(\bz).
\end{gathered}
\end{equation}

Unlike Hopf bifurcation, however, there is no guarantee that the linearisation \eqref{toruspflin} is non-singular since
\[\left[\begin{pmatrix}
		0&-\omz\\ \omz&0
	\end{pmatrix}
	-\mI\pdt-\omz\mI\pdp\right]\bzhtp=\zz\]
may have the solution $\bzhtp\equiv \left(z_1\tp, z_2\tp\right)^T$ with
\[z_1\tp+\ii z_2\tp=\ee^{\ii(\ell\theta+m\phi)}\quad\ell,m\in\Z.\]
This occurs if $\omz(1-m)=\ell$,
and so in particular for $(\ell,m)=(0,1)$; but this is the same as for Hopf bifurcation and again compensated for by the scalar unknowns $\lam,\om$
and the scalar conditions $\bgam$. Now, however, there is a difficulty whenever $\omz$ is rational, i.e. the resonance situation
\begin{equation}\label{pfreson}
	\omz=\frac{\ell}{1-m}.
\end{equation}

One theoretical answer to this problem is to assume that $\omz$ is not only irrational, but also satisfies a Diophantine condition implying that
it is badly approximated by rationals; i.e. $(l,m)$ must be large in order to approximately satisfy \eqref{pfreson}.
This is the approach used in KAM theory \cite{Moser66}, but here we can make a pair of simpler assumptions.
\begin{assumption}\label{pfassum1}No strong resonance, i.e. 
	\[\boxed{\omz\notin\{\tfrac13,\tfrac14\}.}\]
\end{assumption}
\noindent (Here we must remember that our form of Floquet theory in \S\ref{CompFT} enforces the bound $0<\omz<\half$.)
This assumption is required because $\omz$ being rational is also a necessary condition for \ul{subharmonic}
bifurcation of \eqref{pfde} to occur \cite{IJ90}. For rational $\omz$ with denominator $\ge5$, the torus bifurcation is generic; while for $\omz=\frac13$,
the subharmonic bifurcation is generic. (For $\omz=\frac14$, the relative size of certain parameters determines whether torus or subharmonic bifurcation
occurs \cite{IJ90,Wan78r}, but for simplicity we omit this case.)
\begin{theorem}\label{TheoE3pf}
	Under Assumption~\ref{pfassum1}, we can expand \eqref{CRformpf} in powers of $\ep$ and construct an asymptotic solution
	\[\lam^E(\ep),\quad\om^E(\ep)\quad\text{and}\quad\bz^E\tpe\]
	up to and including the $\ep^2$ term, i.e. 
	\begin{equation}\label{E3solpf}
		\begin{gathered}
			\lam^E(\ep)\equiv\lamz+\ep^2\lam_2^E,\quad\om^E(\ep)\equiv\omz+\ep^2\om_2^E,\\
			\bz^E\tpe\equiv\bazp+\ep\bz_2^E\tp+\ep^2\bz_3^E\tp;
		\end{gathered}
	\end{equation}
	where $\bz_2^E$ only depends on the Fourier $\phi$-modes $0$ and $2$ and $\bz_3^E$ only depends on the Fourier $\phi$-modes
	$1$ and $3$. The amplitude and phase conditions force
	\[\left\langle\left\langle\bz_3^E\tp,\bazp\right\rangle\right\rangle=0=\left\langle\left\langle\bz_3^E\tp,\dbazp\right\rangle\right\rangle.\]
	\eqref{E3solpf} can also be expressed in terms of Fourier $\phi$-modes, i.e.
	\begin{equation}\label{E3solpfF}
		\bz^E\tpe\equiv\bazp+\ba_0^E\tep+\sum_{m=1}^3\ba_m^E\tep\cos{m\phi}+\bb_m^E\tep\sin{m\phi},
	\end{equation}
	where $\ba_0^E\tep$, $\ba_2^E\tep$, $\bb_2^E\tep$ are $\ep$-terms and $\ba_1^E\tep$, $\bb_1^E\tep$, $\ba_3^E\tep$, $\bb_3^E\tep$ are $\ep^2$-terms.
	Again, the amplitude and phase conditions force
	\begin{equation}\label{pfap}
		\big\langle\ba_1^E\tep,\bef\big\rangle+\big\langle\bb_1^E\tep,\bet\big\rangle=0
			=\big\langle\ba_1^E\tep,\bet\big\rangle-\big\langle\bb_1^E\tep,\bef\big\rangle.
	\end{equation}
\end{theorem}
Assumption~\ref{pfassum1} is also sufficient to approximately solve \eqref{CRformpf} with $M=3$ Fourier $\phi$-modes; i.e.
\begin{equation}\label{CRpfF3}
\begin{gathered}
	\mS_{\infty,3}\bG\Big(\bz_3\tp,\lam,\theta;\ep\Big)-\pdt\bz_3\tp-\om\pdp\bz_3\tp=\zz\\
	\bgam(\bz_3)-\beh_1=\zz,
\end{gathered}
\end{equation}
where the operator
$\mS_{\infty,3}:{L^2\mapsto L^2}$ performs the Fourier \ul{$\phi$-mode} truncation.
\begin{theorem}\label{TheoF3pf}
	For $\lvert\ep\rvert$ sufficiently small, \eqref{CRpfF3} has a locally unique solution
	\begin{equation}\label{F3solpf}
		\begin{gathered}
			\lam^F(\ep),\;\om^F(\ep)\;\text{and}\\
			\bz^F\tpe\equiv\bazp+\ba_0^F\tep+\sum_{m=1}^3\ba_m^F\tep\cos{m\phi}+\bb_m^F\tep\sin{m\phi}.
		\end{gathered}
	\end{equation}
	As in \eqref{pfap}, the amplitude and phase conditions force
	\begin{equation}\label{Fourpfap}
		\big\langle\ba_1^F\tep,\bef\big\rangle+\big\langle\bb_1^F\tep,\bet\big\rangle=0
			=\big\langle\ba_1^F\tep,\bet\big\rangle-\big\langle\bb_1^F\tep,\bef\big\rangle.
	\end{equation}
	Comparing \eqref{E3solpfF} and \eqref{F3solpf}, as in Corollary~\ref{HopfF3err}, gives the errors
	\begin{align*}
		&&\left\lvert\lam^E(\ep)-\lam^F(\ep)\right\rvert&=O(\ep^4),&\left\lvert\om^E(\ep)-\om^F(\ep)\right\rvert&=O(\ep^4),&\\
		m&=0,2&\left\lVert\ba_m^E\tep-\ba_m^F\tep\right\rVert&=O(\ep^3),&\left\lVert\bb_m^E\tep-\bb_m^F\tep\right\rVert&=O(\ep^3)&m&=2,\\
		m&=1,3&\left\lVert\ba_m^E\tep-\ba_m^F\tep\right\rVert&=O(\ep^4),&\left\lVert\bb_m^E\tep-\bb_m^F\tep\right\rVert&=O(\ep^4)&m&=1,3.
	\end{align*}
\end{theorem}

We can now state our second condition, which may be expressed in several equivalent forms.
\begin{assumption}\label{pfassum2}Nonzero real Lyapunov coefficient, i.e.
	\[\boxed{\lam_2^E\equiv\thalf\left[\lam^F\right]''(0)\neq0.}\]
\end{assumption}
\noindent Since $\lam^E(\ep)-\lamz$ and $\lam^F(\ep)-\lamz$ have no $O(\ep)$ term, Assumption~\ref{pfassum2} forces $\lam^E(\ep)$ and $\lam^F(\ep)$ to
move away from the critical value $\lamz$ for small $\ep\neq0$. Together with Assumption~\ref{pfassum1}, it also shows that $\aR(\lam^E(\ep))$ and
$\aR(\lam^F(\ep))$ move away from zero for small $\ep\neq0$ and therefore permits merely
the no \ul{strong} resonance condition in Assumption~\ref{pfassum1}. (This pair of assumptions has its analogue for Hamiltonian systems \cite{Moser01}.)
We shall see later, in \eqref{exactLpf} and Theorem~\ref{Theoerrpf},
that Assumption~\ref{pfassum2} is equivalent to a real Lyapunov coefficient being nonzero.

\subsection{Normal form and its Fourier approximation}\label{approxnfpf}
In order to cope with possible weak resonances, we need to reduce our equations to an approximate normal form.
Our algorithms in \S\ref{Hopfbif} for the existence, uniqueness and Fourier approximation of periodic orbits created at a Hopf bifurcation point
required neither reduction to the centre manifold nor transformation to normal form: for Neimark--Sacker bifurcation,
however, these two procedures have to be implemented approximately and in this subsection we follow the strategy in \S\ref{approxnf}.

Our aim is to simplify the key equation \eqref{CRformpf}, i.e.
\begin{equation}\label{cmnfdepf}
	\begin{gathered}
		\bG\Big(\bztp,\lam,\theta;\ep\Big)-\left[\pdt+\om\pdp\right]\bztp=\zz\\
		\bgam(\bz)-\beh_1=\zz.
	\end{gathered}
\end{equation}
By again introducing
\[\bz\equiv\left(\bzh,\wtz\right)^T,\quad\text{with $\bzh\in\Rt$ and $\wtz\in\Rnt$,}\]
we can write
\[\bG\Big(\bz,\lam,\theta;\ep\Big)\equiv
	\begin{bmatrix}
		\bGh\Big(\bzh,\wtz,\lam,\theta;\ep\Big)\\
		\wbG\Big(\bzh,\wtz,\lam,\theta;\ep\Big)
	\end{bmatrix}\qquad
	\begin{aligned}
		\bGh:{\Rt\times\Rnt\times\R\times\So\times\R}&\mapsto\Rt\\
		\wbG:{\Rt\times\Rnt\times\R\times\So\times\R}&\mapsto\Rnt.
	\end{aligned}\]
We then construct
\[\bhh:{\Rt\times\So\times\R}\mapsto\Rt\quad\text{and}\quad\wbh:{\Rt\times\So\times\R}\mapsto\Rnt,\]
where $\wbh$ is a homogeneous quadratic polynomial with $\tl$-dependent coefficients and
$\bhh$ is the sum of homogeneous quadratic and cubic polynomials with $\tl$-dependent coefficients, to define the near-identity transformations
\begin{subequations}\label{hhtdefpf}
\begin{equation}
\begin{aligned}
	\bzh&=\byh+\rep\bhh(\ep\byh;\theta,\lam)\\
	&=\byh+\ep\left\{\left[\hy_1^2+\hy_2^2\right]\bah_0\tl+\left[\hy_1^2-\hy_2^2\right]\bah_2\tl+2\hy_1\hy_2\bbh_2\tl\right\}\\
	&\qquad+\ep^2\left\{\hy_1\left[\hy_1^2+\hy_2^2\right]\bah_1\tl+\hy_2\left[\hy_1^2+\hy_2^2\right]\bbh_1\tl\right.\\
		&\qquad\qquad\left.+\hy_1\left[\hy_1^2-3\hy_2^2\right]\bah_3\tl+\hy_2\left[3\hy_1^2-\hy_2^2\right]\bbh_3\tl\right\}\\
\end{aligned}
\end{equation}
and
\begin{equation}
\begin{aligned}
	\wtz&=\wty+\rep\wbh(\ep\byh;\theta,\lam)\\
	&=\wty+\ep\left\{\left[\hy_1^2+\hy_2^2\right]\wba_0\tl+\left[\hy_1^2-\hy_2^2\right]\wba_2\tl+2\hy_1\hy_2\wbb_2\tl\right\}.
\end{aligned}
\end{equation}
\end{subequations}
As in \eqref{hhtdef}, we must also have the restrictions 
\begin{equation}
	\big\langle\bah_1\tl,\beh_1\big\rangle+\big\langle\bbh_1\tl,\beh_2\big\rangle=0
		=\big\langle\bah_1\tl,\beh_2\big\rangle-\big\langle\bbh_1\tl,\beh_1\big\rangle
\end{equation}
in the definition of $\bhh$.
Under these near-identity transformations, \eqref{cmnfdepf} becomes
\begin{subequations}\label{csubdepf}
\begin{gather}
	\bGhd\Big(\byhtp,\wtytp,\lam,\theta;\ep\Big)-\left[\pdt+\om\pdp\right]\byhtp=\zz\label{csubdepfa}\\
	\wbGd\Big(\byhtp,\wtytp,\lam,\theta;\ep\Big)-\left[\pdt+\om\pdp\right]\wtytp=\zz\label{csubdepfb}\\
	\bgam(\byh)-\beh_1=\zz:\label{csubdepfc}
\end{gather}
\end{subequations}
the two mappings
\[\bGhd:{\Rt\times\Rnt\times\R\times\So\times\R}\mapsto\Rt\quad\text{and}\quad\wbGd:{\Rt\times\Rnt\times\R\times\So\times\R}\mapsto\Rnt\]
capable of being expanded, like \eqref{Gexpandpf}, in the form
\begin{align*}
	\bGhd\Big(\byh,\wty,\lam,\theta;\ep\Big)&=\mEhl\byh+\sum_{p\ge2}\ep^{p-1}\bGhd_p\Big(\byh,\wty;\lam,\theta\Big)
		\qquad\bGhd_p:\Rt\times\Rnt\times\R\times\So\mapsto\Rt\\
	\wbGd\Big(\byh,\wty,\lam,\theta;\ep\Big)&=\wEl\wty+\sum_{p\ge2}\ep^{p-1}\wbGd_p\Big(\byh,\wty;\lam,\theta\Big)
		\qquad\wbGd_p:\Rt\times\Rnt\times\R\times\So\mapsto\Rnt;
\end{align*}
where the components of $\bGhd_p$ and $\wbGd_p$ are homogeneous polynomials of degree $p$ in the components of $\byh$ and $\wty$, with coefficients
depending on $\lam$ and $\theta$. Now we choose $\bhh$ and $\wbh$ so that the lower terms in $\bGhd$ and $\wbGd$ may be simplified in the following way:
\begin{itemize}
	\item $\wbh$ forces the coefficients of the quadratic terms for $\byh$ in $\wbGd_2$ to be zero;
	\item $\bhh$ forces the coefficients of the quadratic terms for $\byh$ in $\bGhd_2$ to be zero
		and the coefficients of the cubic terms for $\byh$ in $\bGhd_3$ to take the form
		\begin{equation}\label{Lyappf}
			\left[\hy_1^2+\hy_2^2\right]\mBhl\byh,\qquad\text{where}\quad
			\mBhl\equiv\left(\begin{array}{rr} \bRl&-\bIl\\ \bIl&\bRl \end{array}\right)
		\end{equation}
		and we again call the elements of this matrix Lyapunov coefficients.
\end{itemize}
(I.e. after transformation, only a multiple of the resonant cubic terms \eqref{cubres} remains.)

After this simplification, and under Assumption~\ref{pfassum1}, if we now insert
\begin{subequations}\label{csubsolpf}
\begin{gather}
	\byh\tp=\bahzp\;\left(\equiv(\cos{\phi},\,\sin{\phi})^T\right)\quad\text{and}\quad\wty\tp=\zz\label{csubsolpfa}\\
	\lam=\lamz-\ep^2\frac{\bRlz}{\daRlz}\quad\text{and}\quad\om=\omz+\ep^2\frac{\daRlz\bIlz-\daIlz\bRlz}{\daRlz}\label{csubsolpfb}
\end{gather}
\end{subequations}
into the left-hand side of \eqref{csubdepf}, we can easily see that the remainder is $O(\ep^3)$ for \eqref{csubdepfa}, $O(\ep^2)$ for \eqref{csubdepfb} and
zero for \eqref{csubdepfc}. Consequently, by transforming \eqref{csubsolpfa} back through \eqref{hhtdefpf}, i.e.
\begin{equation}\label{csubsolpftrans}
	\bzh\tp=\bahzp+\rep\bhh(\ep\bahzp;\theta,\lam)\quad\text{and}\quad\wtz\tp=\rep\wbh(\ep\bahzp;\theta,\lam),
\end{equation}
we obtain an asymptotic solution for \eqref{cmnfdepf}. Since Theorem~\ref{TheoE3pf} already displays such a solution, i.e. $\lam^E(\ep)$, $\om^E(\ep)$ and
\[\bz^E\tpe\equiv\begin{pmatrix}\bzh^E\tpe\\ \wtz^E\tpe\end{pmatrix},\]
this must match with \eqref{csubsolpfb} and \eqref{csubsolpftrans}.
Thus we obtain
\begin{equation}\label{cmnfsolpfa}
\begin{aligned}
	\lam^E(\ep)&=\lamz-\ep^2\frac{\bRlz}{\daRlz}\\
	\om^E(\ep)&=\omz+\ep^2\frac{\daRlz\bIlz-\daIlz\bRlz}{\daRlz}
\end{aligned}
\end{equation}
and, through \eqref{csubsolpftrans},
\begin{subequations}\label{cmnfsolpfb}
\begin{align}
	\bzh^E\tpe&=\bahzp+\ep\left\{\bah_0\tlz+\bah_2\tlz\cos{2\phi}+\bbh_2\tlz\sin{2\phi}\right\}\notag\\
		&\qquad+\ep^2\left\{\bah_1\tlz\cos{\phi}+\bbh_1\tlz\sin{\phi}\right.\notag\\
		&\qquad\qquad\left.+\bah_3\tlz\cos{3\phi}+\bbh_3\tlz\sin{3\phi}\right\}\\
	\wtz^E\tpe&=\ep\left\{\wba_0\tlz+\wba_2\tlz\cos{2\phi}+\wbb_2\tlz\sin{2\phi}\right\}+O(\ep^2).
\end{align}
\end{subequations}
Finally, by comparing \eqref{cmnfsolpfb} and \eqref{E3solpf}, we see that the coefficients of $\bhh(.;\theta,\lamz)$ and $\wbh(.;\theta,\lamz)$ are given
\ul{exactly} by the coefficients of the Fourier $\phi$-modes in the $\bzh_2^E\tp$ and $\bzh_3^E\tp$ terms of $\bzh^E\tpe$
and the $\wtz_2^E\tp$ term in $\wtz^E\tpe$ for \eqref{E3solpf}.
Moreover, by comparing \eqref{cmnfsolpfa} and \eqref{E3solpf}, the Lyapunov coefficients $\bRlz$ and $\bIlz$ in \eqref{Lyappf} are given \ul{exactly} by
\begin{equation}\label{exactLpf}
	\bRlz=-\daRlz\lam_2^E\quad\text{and}\quad\bIlz=\om_2^E-\daIlz\lam_2^E,
\end{equation}
and now we see that Assumption~\ref{pfassum2} is equivalent to $\bRlz\neq0$.

To calculate the expansion in \eqref{E3solpf}, however, requires (through $\bG$) explicit knowledge of the second and third derivatives of $\bF$,
so it is practically much more convenient to approximate not only the coefficients of $\bhh(.;\theta,\lamz)$ and $\wbh(.;\theta,\lamz)$ but also the
Lyapunov coefficients $\bRlz$ and $\bIlz$ by using instead the \ul{$M=3$} Fourier $\phi$-mode approximation in \eqref{F3solpf},
i.e. $\lam^F(\ep)$, $\om^F(\ep)$ and
\[\bz^F\tpe\equiv\left\{
\begin{aligned}
	\bzh^F\tpe&\equiv\bah_0^F(\theta,\ep)+\sum_{m=1}^3\bah_m^F(\theta,\ep)\cos{m\phi}+\bbh_m^F(\theta,\ep)\sin{m\phi}\\
	\wtz^F\tpe&\equiv\wba_0^F(\theta,\ep)+\sum_{m=1}^3\wba_m^F(\theta,\ep)\cos{m\phi}+\wbb_m^F(\theta,\ep)\sin{m\phi}.
\end{aligned}
\right.\]
\begin{theorem}\label{Theoerrpf}
	Using the asymptotic error results in Theorem~\ref{TheoF3pf} on page~\pageref{TheoF3pf}, our practical approximate formulae become
	\begin{align*}
		\bRlz&=-\daRlz\frac{\lam^F(\ep)-\lamz}{\ep^2}+O(\ep^2),\\
		\bIlz&=\frac{\om^F(\ep)-\omz}{\ep^2}-\daIlz\frac{\lam^F(\ep)-\lamz}{\ep^2}+O(\ep^2),
	\end{align*} 
	and
	{\renewcommand{\arraystretch}{1.5}
	\[\begin{array}{|l||l|l|}\hline
		m=0,2&\bah_m\tlz=\rep\bah_m^F(\theta,\ep)+O(\ep^2)&\wba_m\tlz=\rep\wba_m^F(\theta,\ep)+O(\ep^2)\\ \hline
		m=2&\bbh_m\tlz=\rep\bbh_m^F(\theta,\ep)+O(\ep^2)&\wbb_m\tlz=\rep\wbb_m^F(\theta,\ep)+O(\ep^2)\\ \hline
		m=1,3&\bah_m\tlz=\rept\bah_m^F(\theta,\ep)+O(\ep^2)&\bbh_m\tlz=\rept\bbh_m^F(\theta,\ep)+O(\ep^2)\\ \hline
	\end{array}\]}
\end{theorem}

We conclude by emphasizing how the $M=3$ Fourier $\phi$-mode approximation plays the same practical role for Neimark--Sacker bifurcation as that
described in the final paragraph of \S\ref{approxnf}.

\subsection{Numerical results}\label{numrespf}
As a numerical example, we use the forced van der Pol equation \cite{GH83,KO00,Reich00}, which may be written in the form \eqref{pfde} as
\begin{equation}\label{fvdpa}
	\dot{x}_1=x_2+\sigma x_1\left(1-x_1^2/3\right),\quad\dot{x}_2=-x_1+\lam\cos{\nu t}.
\end{equation}
Here $\sigma\ge0$ and $0<\nu<1$ are regarded as fixed parameters and $\lam$, as usual, is our continuation parameter: in the form \eqref{pfdea},
\eqref{fvdpa} becomes
\begin{equation}\label{fvdpb}
        \dot{v}_1=\frac1{\nu}\left\{v_2+\sigma v_1\left(1-v_1^2/3\right)\right\},\quad\dot{v}_2=\frac1{\nu}\left\{-v_1+\lam\cos{\theta}\right\}.
\end{equation}
For $\sigma=0$, it is interesting that \eqref{fvdpb} has the periodic orbit and Floquet variables
\begin{equation}\label{fvdpsv}
	\bvt\equiv\frac{\lam}{1-\nu^2}\left(\begin{array}{r}\cos{\theta}\\-\nu\sin{\theta}\end{array}\right),\quad
	\mPt\equiv\begin{pmatrix}1&0\\0&1\end{pmatrix},\quad
	\mE\equiv\frac1{\nu}\left(\begin{array}{rr}0&1\\-1&0\end{array}\right);
\end{equation}
which is useful as a starting value for continuation. (Note that the eigenvalues of $\mE$ in \eqref{fvdpsv} are purely imaginary; and in fact there
is a ``degenerate'' Neimark--Sacker bifurcation here, with respect to the parameter $\sigma$, for which the invariant tori formulae, all at $\sigma=0$,
may be written down exactly. This is of no interest to us.) Having computed a periodic orbit at the value of $\sigma$ we are interested in, we can then
fix $\sigma$ and continue in $\lam$, looking for Neimark--Sacker bifurcation points. We use the techniques described in \cite{Moore05} and,
because of the form of the forcing, the periodic orbits have the symmetry
\[\bv(\theta+\pi)=-\bvt\qquad\forall\theta\in\R;\]
which has the important practical simplification that $\bvt$ need only be approximated by \ul{odd} Fourier modes.
\begin{figure}[t]\begin{center}
	\begin{minipage}{\textwidth}
		\[\begin{array}{|c||c|c|c|c|c|c|}\hline
			\nu&\lamz&\omz&\daRlz&\daIlz&\bRlz&\bIlz\\ \hline
			0.86  &0.4349&  0.1092&  -3.3492&  -1.9037& -0.1154&  -0.1607\\
			0.85  &0.4536&  0.1271&  -2.7032&  -1.3431& -0.0906&  -0.1177\\
			0.84  &0.4726&  0.1444&  -2.2858&  -1.0166& -0.0759&  -0.0933\\
			0.83  &0.4917&  0.1614&  -1.9948&  -0.8072& -0.0663&  -0.0779\\
			0.82  &0.5110&  0.1784&  -1.7808&  -0.6634& -0.0597&  -0.0673\\
			0.81  &0.5304&  0.1953&  -1.6173&  -0.5598& -0.0549&  -0.0597\\
			0.80  &0.5498&  0.2124&  -1.4886&  -0.4822& -0.0513&  -0.0540\\
			0.79  &0.5692&  0.2296&  -1.3849&  -0.4223& -0.0486&  -0.0496\\
			0.78  &0.5886&  0.2471&  -1.2999&  -0.3749& -0.0465&  -0.0460\\
			0.77  &0.6079&  0.2648&  -1.2290&  -0.3367& -0.0449&  -0.0432\\
			0.76  &0.6271&  0.2829&  -1.1693&  -0.3054& -0.0436&  -0.0408\\
			0.75  &0.6462&  0.3013&  -1.1185&  -0.2793& -0.0425&  -0.0389\\
			0.74  &0.6651&  0.3200&  -1.0749&  -0.2573& -0.0417&  -0.0372\\
			0.73  &0.6840&  0.3392&  -1.0372&  -0.2386& -0.0411&  -0.0358\\ \hline
		\end{array}\]
	\end{minipage}\end{center}
	\caption{\label{vdpoutput}Neimark--Sacker bifurcation points for the forced van der Pol equation}
\end{figure}
This symmetry is inherited by the Floquet decomposition in \S\ref{CompFT}, so that (if we use the strategy in \cite{Moore05} to limit the size of the
imaginary part of the Floquet exponents) either
\[\mP(\theta+\pi)=\pm\mPt\quad\text{and}\quad\bw(\theta+\pi)=\mp\bwt\qquad\forall\theta\in\R.\]
In Figure~\ref{vdpoutput} we display $\lamz$ for Neimark--Sacker bifurcation points at different $\nu$ values but with $\sigma=4$, and this
may be compared with Figure 13 in \cite{Reich00}. (A simple secant iteration was used to locate the periodic orbits with purely imaginary Floquet
exponents, so we are not using a sophisticated method to detect Neimark--Sacker bifurcation points.) We want to show how some of the important scalars
associated with the bifurcation vary with $\nu$ in this example; and so
we display the $\omz$, $\daRlz$, $\daIlz$, $\bRlz$ and $\bIlz$ values at these bifurcation points, the latter pair being approximated
as in Theorem~\ref{Theoerrpf} with $\ep=0.005$. (Note that we have jumped across two points of strong resonance, where $\omz=\frac14$ and $\frac13$.)
For these calculations we used $M=24$ Fourier $\theta$-modes, which reduced the size of the Fourier coefficients to $\approx 10^{-14}$.

\subsection{Higher-order Fourier approximation of tori}\label{hoFpf}
In order to compute higher-order approximations for our invariant tori, we must employ a more suitable parametrisation than \eqref{toruspf1}.
Thus we use the normal bundle of the approximate torus
\begin{equation}\label{toruspf2}
	\bvztl+\ep\mPtl\bazp
\end{equation}
and, in \eqref{csubdepf},
\begin{itemize}
	\item replace $\byhtp$ with $\left[1+\rhotp\right]\bahzp$ for unknown $\rho:\So\times\So\mapsto\R$,
	\item allow $\om:\So\times\So\mapsto\R$ to be an unknown \ul{function}.
\end{itemize}
This links up with the invariance condition used in \cite{Moore96} for continuation of tori, and corresponds to using polar co-ordinates in the critical
$2$-dimensional subspace. In \eqref{csubdepfc} there is now no need for a scalar phase condition, and the scalar amplitude equation simplifies to a
zero-mean condition for $\rho$, i.e.
\begin{equation}\label{zerompf}
	\left\langle\left\langle\rhotp,1\right\rangle\right\rangle=0.
\end{equation}
Thus our equations for $\rho$ and $\om$ in \eqref{csubdepfa} decouple to become
\begin{subequations}\label{centrepf}
\begin{gather}
	\begin{split}
		\bahzp\cdot\bGhd\Big(\left[1+\rhotp\right]\bahzp,&\wtytp,\lam,\theta;\ep\Big)\\
		&-\left[\pdt+\omtp\pdp\right]\rhotp=0\label{centrepfa}
	\end{split}\\
\intertext{and}
	\frac1{1+\rhotp}\badhzp\cdot\bGhd\Big(\left[1+\rhotp\right]\bahzp,\wtytp,\lam,\theta;\ep\Big)-\omtp=0,\label{centrepfb}
\end{gather}
\end{subequations}
while the hyperbolic equations in \eqref{csubdepfb} remain
\begin{equation}\label{hyperpf}
	\wbGd\Big(\left[1+\rhotp\right]\bahzp,\wtytp,\lam,\theta;\ep\Big)-\left[\pdt+\omtp\pdp\right]\wtytp=\zz.
\end{equation}
The crucial leading terms in \eqref{centrepf} are
\begin{subequations}\label{ltermpf}
\begin{equation}
\begin{split}
	\bahzp\cdot\bGhd\Big(&\left[1+\rhotp\right]\bahzp,\wtytp,\lam,\theta;\ep\Big)\\
	&=\left[1+\rhotp\right]\aRl+\ep^2\bRl\left[1+\rhotp\right]^3+O(\ep^3)
\end{split}
\end{equation}
and
\begin{equation}\label{ltermpfb}
\begin{split}
	\frac1{1+\rhotp}\badhzp\cdot\bGhd\Big(&\left[1+\rhotp\right]\bahzp,\wtytp,\lam,\theta;\ep\Big)\\
	&=\aIl+\ep^2\bIl\left[1+\rhotp\right]^2+O(\ep^3).
\end{split}
\end{equation}
\end{subequations}

Consequently, if we use \eqref{centrepfb} to \ul{define} $\omtp$ in terms of $\lam$, $\rhotp$ and $\wtytp$ for the rest of this subsection,
we finally have to prove that the system of equations
\begin{subequations}\label{fineqpf}
\begin{gather}
	\begin{split}
		\bahzp\cdot\bGhd\Big(\left[1+\rhotp\right]\bahzp,&\wtytp,\lam,\theta;\ep\Big)\\
		&-\left[\pdt+\omtp\pdp\right]\rhotp=0
	\end{split}\\
	\begin{split}
		\wbGd\Big(\left[1+\rhotp\right]\bahzp,&\wtytp,\lam,\theta;\ep\Big)\\
		&-\left[\pdt+\omtp\pdp\right]\wtytp=\zz
	\end{split}\\
	\left\langle\left\langle\rhotp,1\right\rangle\right\rangle=0
\end{gather}
\end{subequations}
has a locally unique solution $\left(\lam,\rho,\wty\right)$ for $\lvert\ep\rvert$ sufficiently small.
This is achieved in \cite{Sack64,Sack65} through the iteration
\begin{subequations}\label{iterpf}
\begin{gather}
	\daRz\delamko+\left\{2\ep^2\bRlz-\left[\pdt+\om^{(k)}\tp\pdp\right]\right\}\rho^{(k+1)}\tp=\hrktp\label{iterpfa}\\
	\left\{\wEz-\left[\pdp+\om^{(k)}\tp\pdp\right]\right\}\wty^{(k+1)}\tp=\wbrktp\\
	\left\langle\left\langle\rho^{(k+1)}\tp,1\right\rangle\right\rangle=0,
\end{gather}
\end{subequations}
where $\delamko\equiv\lamko-\lamk$ and
\begin{align*}
	\om^{(k)}\tp&\equiv\frac1{1+\rho^{(k)}\tp}\badhzp\cdot\bGhd\Big(\left[1+\rho^{(k)}\tp\right]\bahzp,\wty^{(k)}\tp,\lamk,\theta;\ep\Big)\\
	\hrktp&\equiv 2\ep^2\bRlz\rho^{(k)}\tp\\
		&\qquad\qquad-\bahzp\cdot\bGhd\Big(\left[1+\rho^{(k)}\tp\right]\bahzp,\wty^{(k)}\tp,\lamk,\theta;\ep\Big)\\
	\wbrktp&\equiv\wEz\wty^{(k)}\tp-\wbGd\Big(\left[1+\rho^{(k)}\tp\right]\bahzp,\wty^{(k)}\tp,\lamk,\theta;\ep\Big),
\end{align*}
with starting values
\[\lams=\lamz,\;\rho^{(0)}\tp=0,\;\wty^{(0)}\tp=\zz.\]
The key idea behind showing that these iterates remain bounded
and then converge is to integrate \eqref{iterpfa} against $\rho^{(k+1)}\tp$, so that the l.h.s. becomes
\begin{multline}\label{definitepf}
	\left\langle\left\langle\daRz\delamko+\left\{2\ep^2\bRlz-\left[\pdt+\om^{(k)}\tp\pdp\right]\right\}\rho^{(k+1)}\tp,
		\rho^{(k+1)}\tp\right\rangle\right\rangle\\
	=\left\langle\left\langle\left\{2\ep^2\bRlz+\pdp\om^{(k)}\tp\right\}\rho^{(k+1)}\tp,\rho^{(k+1)}\tp\right\rangle\right\rangle.
\end{multline}
Since \eqref{ltermpfb} shows that the leading \ul{non-constant} term in $\om^{(k)}\tp$ is $O(\ep^3)$, Assumption~\ref{pfassum2} ensures that
\eqref{definitepf} is a definite quadratic term for $\lvert\ep\rvert$ sufficiently small, and this is sufficient for \cite{Sack64} to prove the
following result.
\begin{theorem}\label{Sackpf}
	Suppose $\bF$ in \eqref{pfdea} has $r\ge5$ continuous derivatives for $(\lam,\bv)$ in a neighbourhood of $(\lamz,\bvzt)$. Then $\exists\,{\ep_r>0}$
	such that for $\lvert\ep\rvert<\ep_r$ \eqref{fineqpf} has a locally unique solution
	\[\lamz(\ep),\;\rhoz\tpe,\;\wtyz\tpe\]
	with $\rhoz(.,.;\ep),\;\wtyz(.,.;\ep)$ having $(r-1)$ Lipchitz continuous derivatives. This means that both $\byh$ and $\bzh,\;\wtz$
	(through \eqref{hhtdefpf}) have this degree of smoothness, and so, through \eqref{toruspf1}, do the invariant tori as manifolds.
\end{theorem}

\noindent The subtlety of Theorem~\ref{Sackpf} is that, in general, $\ep_r\to 0$ as $r\to\infty$; in particular, one cannot expect the tori to be analytic
when $\bF$ is analytic.

In practice we seek an approximate solution of \eqref{fineqpf} in the form
\[\begin{aligned}
	\rhoLM\tp&\equiv\sum_{\ell,m}\rho_{\ell,m}\ee^{\ii(\ell\theta+m\phi)}\\
	\wtyLM\tp&\equiv\sum_{\ell,m}\wty_{\ell,m}\ee^{\ii(\ell\theta+m\phi)}
\end{aligned}\qquad\lvert\ell\rvert\le L,\;\lvert m\rvert\le M\]
with $\rho_{-\ell,-m},\;\wty_{-\ell,-m}$ the conjugates of $\rho_{\ell,m},\;\wty_{\ell,m}$. These functions must satisfy
\begin{subequations}\label{fineqFpf}
\begin{gather}
\begin{split}
	\left\langle\left\langle\bahzp\cdot\bGhd\Big(\right.\right.&\left[1+\rhoLM\tp\right]\bahzp,\wtyLM\tp,\lam,\theta;\ep\Big)\\
	&\left.\left.-\left[\pdt+\omtp\pdp\right]\rhoLM\tp,\ee^{\ii(\ell\theta+m\phi)}\right\rangle\right\rangle=0
\end{split}\\
\begin{split}
	\left\langle\left\langle\wbGd\Big(\right.\right.&\left[1+\rhoLM\tp\right]\bahzp,\wtyLM\tp,\lam,\theta;\ep\Big)\\
	&\left.\left.-\left[\pdt+\omtp\pdp\right]\wtyLM\tp,\bx\ee^{\ii(\ell\theta+m\phi)}\right\rangle\right\rangle=\zz
\end{split}\\
	\rho_{0,0}=0,
\end{gather}
\end{subequations}
for $\lvert\ell\rvert\le L$, $\lvert m\rvert\le M$ and $\forall{\bx\in\Rnt}$, with
\[\omtp\equiv\frac1{1+\rhoLM\tp}\badhzp\cdot\bGhd\Big(\left[1+\rhoLM\tp\right]\bahzp,\wtyLM\tp,\lam,\theta;\ep\Big).\]
As shown in \cite{Sam91}, the iteration analogous to \eqref{iterpf} also converges here and gives the following result.
\begin{theorem}\label{Sampf}
	Under the conditions of Theorem~\ref{Sackpf}, \eqref{fineqFpf} has a locally unique solution
	\[\lamLM^F(\ep),\;\rhoLM^F\tpe,\;\wtyLM^F\tpe\]
	satisfying
	\[\max\left\{\begin{array}{c}
		\lvert\lamLM^F(\ep)-\lamz(\ep)\rvert\\
		\lVert\rhoLM^F(.,.;\ep)-\mSLM\rhoz(.,.;\ep)\rVert_{L_2}\\
		\lVert\wtyLM^F(.,.;\ep)-\mSLM\wtyz(.,.;\ep)\rVert_{L_2}
	\end{array}\right\}\le C
	\max\left\{\begin{array}{c}
		\lVert\left(\mI-\mSLM\right)\rhoz(.,.;\ep)\rVert_{H^1}\\
		\lVert\left(\mI-\mSLM\right)\wtyz(.,.;\ep)\rVert_{H^1}
	\end{array}\right\}.\]
\end{theorem}
We comment on the implementation of this algorithm in \S\ref{Conclu}.

\section{Neimark--Sacker bifurcation for autonomous systems}\label{NSbifaut}
We start with our two basic conditions.
\begin{assumption}\label{NSautfirst}
	At $\lam=\lamz$, \eqref{autde} has a periodic orbit $\buz(t)$ of period $\tpTz$, and so, under the change-of-variable
	\[\bvzt\equiv\buz(\theta\Tz)\qquad\bvz:{\So\mapsto\Rn},\]
	$\Tz$ and $\bvz$ satisfy
	\[\Tz\bFvztlz-\ddt\bvzt=\zz.\]
\end{assumption}
\begin{assumption}
	If we apply the Floquet theory in \S\ref{CompFT} to
	\[\mAt\equiv\Tz\mJvztlz\]
	then \eqref{Floquev} becomes
	\[\Tz\mJvztlz\mPzt-\dPzt=\mPzt\mEz,\]
	where $\mEz\in\Rnn$ and $\mPz:\So\mapsto\Rnn$ with $\mPzt$ non-singular $\forall\theta\in\So$, and we have the invariant subspace decomposition
	\[\mEz\equiv\left(\begin{array}{lll}
				\mEhz&\ZZ&\zz\\
				\ZZ&\wEz&\zz\\
				\zz&\zz&0
		\end{array}\right)\quad
		\begin{aligned}
			\mEhz&\in\Rtt\\
			\wEz&\in\Rnthnth\\
		\end{aligned}
	\qquad\text{with}\quad
	\mEhz\equiv\begin{pmatrix}
			0&-\omz\\ \omz&0
			\end{pmatrix}\quad\omz>0\]
	and $\wEz$ having no eigenvalues on the imaginary axis.
\end{assumption}

\noindent The \IFT\ then gives us a locally unique curve of periodic orbits, smoothly parametrised by $\lam$, and satisfying
\begin{equation}\label{autcurve}
	\Tzl\bFvztll-\ddt\bvztl=\zz.
\end{equation}
The Floquet variables in the invariant subspace decomposition can be smoothly continued locally, and so we have
\begin{equation}\label{autinvar}
	\Tzl\mJvztll\mPtl-\dPtl=\mPtl\mEl\qquad\begin{gathered}\mP:{{\So\times\R}\mapsto\Rnn}\\ \mE:{\R\mapsto\Rnn},\end{gathered}
\end{equation}
where $\mPtl$ is non-singular and
\[\mEl\equiv\begin{pmatrix}
		\mEhl&\ZZ&\zz\\
		\ZZ&\wEl&\zz\\
		\zz&\zz&0
	\end{pmatrix}\qquad\begin{aligned}
				\mEh&:\R\mapsto\Rtt\\
				\wE&:\R\mapsto\Rnthnth\\
			\end{aligned}\]
with
\[\mEhl\equiv\left(\begin{array}{rr}
			\aRl&-\aIl\\ \aIl&\aRl
		\end{array}\right)\qquad\begin{aligned}
						\aR&:\R\mapsto\R\\
						\aI&:\R\mapsto\R.
					\end{aligned}\]
Finally, the key transversality condition must also hold.
\begin{assumption}\label{auttrans}Transversal crossing of critical Floquet exponents, i.e.
	\[\boxed{\daRz\equiv\daRlz\neq0.}\]
\end{assumption}

\subsection{Crandall--Rabinowitz formulation}\label{secCRformaut}
To start with, we attempt to mimic our approach in \S\ref{secCRformpf} and seek invariant tori of \eqref{autde} in the form
\begin{equation}\label{torusaut1}
	\bvztl+\ep\mPtl\bztp\qquad\bz:\So\times\So\mapsto\Rn,
\end{equation}
with unknown $\bz$, satisfying
\begin{subequations}\label{invaraut1}
\begin{multline}\label{invaraut1a}
	\Tzl\bF\Big(\bvztl+\ep\mPtl\bztp,\lam\Big)\\
	-\left[1+\ep\eta\right]\pdt\Big[\bvztl+\ep\mPtl\bztp\Big]\\
	-\om\pdp\Big[\bvztl+\ep\mPtl\bztp\Big]=\zz
\end{multline}
for unknown $\om,\eta\in\R$. As in \S\ref{secCRformpf}, we are expressing the invariance of \eqref{torusaut1} by insisting that the vector
field lie in its tangent space; the only difference being the extra unknown $\eta$ now compensating for the zero Floquet exponent.
The parametrisation of the torus in \eqref{invaraut1} is again special, however, since the coefficients ($[1+\ep\eta]$
and $\om$) must be constant: as in \S\ref{secCRformpf}, we will need to generalise this parametrisation later in \S\ref{hoFaut}.
Of course, we also require the scalar amplitude and phase conditions
\begin{equation}\label{autscona}
	\bgam_1(\bz)\equiv
		\begin{pmatrix}
			\left\langle\left\langle\bztp,\bazp\right\rangle\right\rangle\\
			\left\langle\left\langle\bztp,\dbazp\right\rangle\right\rangle
		\end{pmatrix}=\beh_1,
\end{equation}
and the scalar phase condition
\begin{equation}\label{autsconb}
	\gam_2(\bz)\equiv\left\langle\left\langle\bztp,\ben\right\rangle\right\rangle=0.
\end{equation}
\end{subequations}

To attempt to apply the \IFT\ to \eqref{invaraut1}, we must first eliminate the curve of periodic orbits \eqref{autcurve}: thus the
Crandall--Rabinowitz formulation writes
\begin{multline*}
	\bG\Big(\bztp,\lam,\eta,\theta;\ep\Big)\equiv\rep\mPtl^{-1}\bigg\{\Tzl\bF\Big(\bvztl+\ep\mPtl\bztp,\lam\Big)\\
	-\Tzl\bF\Big(\bvztl,\lam\Big)-\ep\left[1+\ep\eta\right]\dPtl\bztp\bigg\}
\end{multline*}
and solves \eqref{invaraut1} in the form
\begin{equation}\label{CRformaut1}
	\zz=\cF(\bztp,\lam,\eta,\om,\theta;\ep)\equiv\\ \left\{
	\begin{gathered}
		\begin{aligned}
			\bG(\bztp,&\lam,\eta,\theta;\ep)-\eta\ben\\
			&-\left[(1+\ep\eta)\pdt+\om\pdp\right]\bztp
		\end{aligned}\\
		\bgam_1(\bz)-\beh_1\\
		\gam_2(\bz).
	\end{gathered}\right.
\end{equation}
Hence, using \eqref{autcurve} and \eqref{autinvar}, we can expand $\bG$ in the form
\begin{equation}\label{Gexpandaut}
	\bG\Big(\bz,\lam,\eta,\theta;\ep\Big)=\mEl\bz+\ep\eta\bG_1\Big(\bz;\lam,\theta\Big)+\sum_{p\ge2}\ep^{p-1}\bG_p\Big(\bz;\lam,\theta\Big),
\end{equation}
where $\bG_p:\Rn\times\R\times\So\mapsto\Rn$ and
the $n$ components of $\bG_p$ are homogeneous polynomials of degree $p$ in the $n$ components of $\bz$ with coefficients depending on $\lam$ and $\theta$.
At $\ep=0$, \eqref{CRformaut1} becomes
\begin{gather*}
	\Big[\mEl-\mI\pdt-\om\mI\pdp\Big]\bztp-\eta\ben=\zz\\
	 \bgam_1(\bz)-\beh_1=\zz\\
	 \gam_2(\bz)=0
\end{gather*}
with solution 
\[\bztp=\bazp,\;\lam=\lamz,\;\om=\omz,\;\eta=0;\]
and the linearisation about this solution is
\begin{equation}\label{torusaut1lin}
\begin{gathered}
	\left[\mEz-\mI\pdt-\omz\mI\pdp\right]\bztp+\lam\dElz\bazp-\om\dbazp-\eta\ben\\
	\bgam_1(\bz)\\
	\gam_2(\bz).
\end{gathered}
\end{equation}

Just as in \S\ref{secCRformpf}, there is no guarantee that the linearisation \eqref{torusaut1lin} is non-singular, and singularity occurs if
\begin{equation}\label{resonaut}
	(m-1)\omz+\ell=0\quad\text{or}\quad m\omz+\ell=0\qquad\text{for some $\ell,m\in\Z$}.
\end{equation}
The first occurs for $(\ell,m)=(0,1)$, but this is compensated for by the scalar unknowns $\lam,\om$ and the scalar conditions $\bgam_1$; and the
second occurs for $(\ell,m)=(0,0)$, but this is compensated for by the scalar unknown $\eta$ and the scalar condition $\gam_2$. If $\omz$ is rational,
however, \eqref{resonaut} will be satisfied by larger integer values of $(\ell,m)$; even if $\omz$ is irrational, they will be satisfied ``arbitrarily''
closely. Thus we must impose the same condition as in in \S\ref{secCRformpf}.
\begin{assumption}\label{autassum1}No strong resonance, i.e.
	\[\boxed{\omz\notin\{\tfrac13,\tfrac14\}.}\]
\end{assumption}
\noindent It then follows that an asymptotic solution for \eqref{CRformaut1} can be constructed. 
\begin{theorem}\label{TheoE3aut}
	Under Assumption~\ref{autassum1}, we can expand \eqref{CRformaut1} in powers of $\ep$ and construct an asymptotic solution
	\[\lam^E(\ep),\quad\om^E(\ep),\quad\eta^E(\ep)\quad\text{and}\quad\bz^E\tpe\]
	up to and including the $\ep^2$ terms, i.e. 
	\begin{equation}\label{E3solaut}
		\begin{gathered}
			\lam^E(\ep)\equiv\lamz+\ep^2\lam_2^E,\quad\om^E(\ep)\equiv\omz+\ep^2\om_2^E,\quad\eta^E(\ep)\equiv\ep\eta_2^E,\\
			\bz^E\tpe\equiv\bazp+\ep\bz_2^E\tp+\ep^2\bz_3^E\tp;
		\end{gathered}
	\end{equation}
	where $\bz_2^E$ only depends on the Fourier $\phi$-modes $0$ and $2$ and $\bz_3^E$ only depends on the Fourier $\phi$-modes
	$1$ and $3$. The amplitude and phase conditions force
	\[\left\langle\left\langle\bz_3^E\tp,\bazp\right\rangle\right\rangle=0=\left\langle\left\langle\bz_3^E\tp,\dbazp\right\rangle\right\rangle
		\quad\text{and}\quad\left\langle\left\langle\bz_2^E\tp,\ben\right\rangle\right\rangle=0.\]
	\eqref{E3solaut} can also be expressed in terms of Fourier $\phi$-modes, i.e.
	\begin{equation}\label{E3solautF}
		\bz^E\tpe\equiv\bazp+\ba_0^E\tep+\sum_{m=1}^3\ba_m^E\tep\cos{m\phi}+\bb_m^E\tep\sin{m\phi},
	\end{equation}
	where $\ba_0^E\tep$, $\ba_2^E\tep$, $\bb_2^E\tep$ are $\ep$-terms and $\ba_1^E\tep$, $\bb_1^E\tep$, $\ba_3^E\tep$, $\bb_3^E\tep$ are $\ep^2$-terms.
	Again, the amplitude and phase conditions force
	\begin{equation}\label{autap}
		\begin{aligned}
			\left\langle\ba_1^E\tep,\bef\right\rangle+\left\langle\bb_1^E\tep,\bet\right\rangle&=0\\
			\left\langle\ba_1^E\tep,\bet\right\rangle-\left\langle\bb_1^E\tep,\bef\right\rangle&=0
		\end{aligned}
		\qquad\text{and}\quad\left\langle\ba_0^E\tep,\ben\right\rangle=0.
	\end{equation}
\end{theorem}
Assumption~\ref{autassum1} is also sufficient to approximately solve \eqref{CRformaut1} with $M=3$ Fourier $\phi$-modes; i.e.
\begin{equation}\label{CRautF3}
\begin{gathered}
	\mS_{\infty,3}\bG\Big(\bz_3\tp,\lam,\eta,\theta;\ep\Big)-\eta\ben-\left[(1+\ep\eta)\pdt+\om\pdp\right]\bz_3\tp=\zz\\
	\bgam_1(\bz_3)-\beh_1=\zz\\
	\gam_2(\bz_3)=0,
\end{gathered}
\end{equation}
where the operator $\mS_{\infty,3}:{L^2\mapsto L^2}$ performs the Fourier \ul{$\phi$-mode} truncation.
\begin{theorem}\label{TheoF3aut}
	For $\lvert\ep\rvert$ sufficiently small, \eqref{CRautF3} has a locally unique solution
	\begin{equation}\label{F3solaut}
		\begin{gathered}
			\lam^F(\ep),\;\om^F(\ep),\;\eta^F(\ep)\;\text{and}\\
			\bz^F\tpe\equiv\bazp+\ba_0^F\tep+\sum_{m=1}^3\ba_m^F\tep\cos{m\phi}+\bb_m^F\tep\sin{m\phi}.
		\end{gathered}
	\end{equation}
	As in \eqref{autap}, the amplitude and phase conditions force
	\begin{equation}\label{Fourautap}
		\begin{aligned}
			\left\langle\ba_1^F\tep,\bef\right\rangle+\left\langle\bb_1^F\tep,\bet\right\rangle&=0\\
			\left\langle\ba_1^F\tep,\bet\right\rangle-\left\langle\bb_1^F\tep,\bef\right\rangle&=0
		\end{aligned}
		\qquad\text{and}\quad\left\langle\ba_0^F\tep,\ben\right\rangle=0.
	\end{equation}
	Comparing \eqref{E3solautF} and \eqref{F3solaut} gives the errors $\left\lvert\eta^E(\ep)-\eta^F(\ep)\right\rvert=O(\ep^3)$,
	\begin{align*}
		&&\left\lvert\lam^E(\ep)-\lam^F(\ep)\right\rvert&=O(\ep^4),&\left\lvert\om^E(\ep)-\om^F(\ep)\right\rvert&=O(\ep^4),&\\
		m&=0,2&\left\lVert\ba_m^E\tep-\ba_m^F\tep\right\rVert&=O(\ep^3),&\left\lVert\bb_m^E\tep-\bb_m^F\tep\right\rVert&=O(\ep^3)&m&=2,\\
		m&=1,3&\left\lVert\ba_m^E\tep-\ba_m^F\tep\right\rVert&=O(\ep^4),&\left\lVert\bb_m^E\tep-\bb_m^F\tep\right\rVert&=O(\ep^4)&m&=1,3.
	\end{align*}
\end{theorem}

We can now state our second condition.
\begin{assumption}\label{autassum2}Nonzero real Lyapunov coefficient, i.e.
	\[\boxed{\lam_2^E\equiv\thalf\left[\lam^F\right]''(0)\neq0,}\]
\end{assumption}
\noindent As in \S\ref{secCRformpf}, this means that $\lam^E(\ep)$ and $\lam^F(\ep)$ move away from the critical value $\lamz$ for small
$\ep\neq0$. \eqref{exactLaut} and Theorem~\ref{Theoerraut} show that this is equivalent to a real Lyapunov coefficient being nonzero.

\subsection{Normal form and its Fourier approximation}\label{approxnfaut}
We follow the strategy in \S\ref{approxnfpf}, and construct the necessary transformations in order to simplify the key equation \eqref{CRformaut1},
i.e.
\begin{equation}\label{cmnfdeaut}
	\begin{gathered}
		\bG\Big(\bztp,\lam,\eta,\theta;\ep\Big)-\eta\ben-\left[(1+\ep\eta)\pdt+\om\pdp\right]\bztp=\zz\\
		\bgam_1(\bz)-\beh_1=\zz\\
		\gam_2(\bz)=0.
	\end{gathered}
\end{equation}
By introducing
\[\bz\equiv\left(\bzh,\wtz,\zo\right)^T,\quad\text{with $\left(\bzh,\wtz,\zo\right)\in{\Rt\times\Rnth\times\R}$,}\]
we can write
\[\bG\Big(\bz,\lam,\eta,\theta;\ep\Big)\equiv
	\begin{bmatrix}
		\bGh\Big(\bzh,\wtz,\zo,\lam,\eta,\theta;\ep\Big)\\
		\wbG\Big(\bzh,\wtz,\zo,\lam,\eta,\theta;\ep\Big)\\
		\Go\Big(\bzh,\wtz,\zo,\lam,\eta,\theta;\ep\Big)
	\end{bmatrix}\quad
	\begin{aligned}
		\bGh:{\Rt\times\Rnth\times\R\times\R\times\R\times\So\times\R}&\mapsto\Rt\\
		\wbG:{\Rt\times\Rnth\times\R\times\R\times\R\times\So\times\R}&\mapsto\Rnth\\
		\Go:{\Rt\times\Rnth\times\R\times\R\times\R\times\So\times\R}&\mapsto\R.
	\end{aligned}\]
We then construct
\[\bhh:{\Rt\times\So\times\R}\mapsto\Rt,\quad\wbh:{\Rt\times\So\times\R}\mapsto\Rnth\quad\text{and}\quad\ho:{\Rt\times\So\times\R}\mapsto\R,\]
where $\ho$ is a homogeneous quadratic polynomial with $\tl$-dependent coefficients and $\bhh$, $\wbh$ are the same as in \S\ref{approxnfpf},
to define the near-identity transformations
\begin{subequations}\label{hhtodefaut}
\begin{align}
	\bzh&=\byh+\rep\bhh(\ep\byh;\theta,\lam)\quad\text{with expansion as in \eqref{hhtdefpf}}\\
	\wtz&=\wty+\rep\wbh(\ep\byh;\theta,\lam)\quad\text{with expansion as in \eqref{hhtdefpf}}\\
	\zo&=\yo+\rep\ho(\ep\byh;\theta,\lam)\notag\\
	&=\yo+\ep\left\{\left[\hy_1^2+\hy_2^2\right]\oa_0\tl+\left[\hy_1^2-\hy_2^2\right]\oa_2\tl+2\hy_1\hy_2\ob_2\tl\right\}.
\end{align}
\end{subequations}
Now we must have the restrictions 
\begin{equation}
	\left\langle\bah_1\tl,\beh_1\right\rangle+\left\langle\bbh_1\tl,\beh_2\right\rangle=0
		=\left\langle\bah_1\tl,\beh_2\right\rangle-\left\langle\bbh_1\tl,\beh_1\right\rangle
\end{equation}
in the definition of $\bhh$, and the restriction
\begin{equation}\label{h0restrict}
	\left\langle\oa_0\tl,1\right\rangle=0
\end{equation}
in the definition of $\ho$.
Finally, to compensate for \eqref{h0restrict}, it is also necessary to include the near-identity transformation
\begin{equation}\label{defkappa}
	\eta=\zeta+\ep\kapl.
\end{equation}

Under all these transformations, \eqref{cmnfdeaut} becomes
\begin{subequations}\label{csubdeaut}
\begin{gather}
	\begin{split}
		\bGhd\Big(\byhtp,&\wtytp,\yotp,\lam,\zeta,\theta;\ep\Big)\\
		&-\left[(1+\ep^2\kapl+\ep\zeta)\pdt+\om\pdp\right]\byhtp=\zz
	\end{split}\label{csubdeauta}\\
	\begin{split}
		\wbGd\Big(\byhtp,&\wtytp,\yotp,\lam,\zeta,\theta;\ep\Big)\\
		&-\left[(1+\ep^2\kapl+\ep\zeta)\pdt+\om\pdp\right]\wtytp=\zz
	\end{split}\label{csubdeautb}\\
	\begin{split}
		\God\Big(\byhtp,&\wtytp,\yotp,\lam,\zeta,\theta;\ep\Big)-\left[\zeta+\ep\kapl\right]\\
		&-\left[(1+\ep^2\kapl+\ep\zeta)\pdt+\om\pdp\right]\yotp=0
	\end{split}\label{csubdeautc}\\
	\bgam_1(\byh)-\beh_1=\zz\label{csubdeautd}\\
	\gam_2(\yo)=0:\label{csubdeaute}
\end{gather}
\end{subequations}
the three mappings
\begin{align*}
	\bGhd&:{\Rt\times\Rnth\times\R\times\R\times\R\times\So\times\R}\mapsto\Rt,\\
	\wbGd&:{\Rt\times\Rnth\times\R\times\R\times\R\times\So\times\R}\mapsto\Rnth\\
	\God&:{\Rt\times\Rnth\times\R\times\R\times\R\times\So\times\R}\mapsto\R
\end{align*}
capable of being expanded, like \eqref{Gexpandaut}, in the form
\begin{align*}
	\bGhd\Big(\byh,\wty,\yo,\lam,\zeta,\theta;\ep\Big)&=\mEhl\byh+\sum_{p\ge2}\ep^{p-1}\bGhd_p\Big(\byh,\wty,\yo;\lam,\theta\Big)
		+\zeta\sum_{p\ge1}\ep^p\bKhd_p\Big(\byh,\wty,\yo;\lam,\theta\Big)\\
	\wbGd\Big(\byh,\wty,\yo,\lam,\zeta,\theta;\ep\Big)&=\wEl\wty+\sum_{p\ge2}\ep^{p-1}\wbGd_p\Big(\byh,\wty,\yo;\lam,\theta\Big)
		+\zeta\sum_{p\ge1}\ep^p\wbKd_p\Big(\byh,\wty,\yo;\lam,\theta\Big)\\
	\God\Big(\byh,\wty,\yo,\lam,\zeta,\theta;\ep\Big)&=\sum_{p\ge2}\ep^{p-1}\God_p\Big(\byh,\wty,\yo;\lam,\theta\Big)
		+\zeta\sum_{p\ge1}\ep^p\Kod_p\Big(\byh,\wty,\yo;\lam,\theta\Big);
\end{align*}
where the components of both $\bGhd_p$, $\wbGd_p$ $\God_p$ and $\bKhd_p$, $\wbKd_p$ $\Kod_p$
are homogeneous polynomials of degree $p$ in the components of $\byh$, $\wty$ and $\yo$, with coefficients depending on $\lam$ and $\theta$.
Now we choose $\bhh$, $\wbh$ and $\ho$ so that the lower terms in $\bGhd$, $\wbGd$ and $\God$ may be simplified in the following way.
\begin{itemize}
	\item $\wbh$ forces the coefficients of the quadratic terms for $\byh$ in $\wbGd_2$ to be zero;
	\item $\ho$, and $\kapl$ in \eqref{defkappa}, force the coefficients of the quadratic terms for $\byh$ in $\God_2$ to take the form
		\[\kapl\left[\hy_1^2+\hy_2^2\right];\]
	\item $\bhh$ forces the coefficients of the quadratic terms for $\byh$ in $\bGhd_2$ to be zero
		and the coefficients of the cubic terms for $\byh$ in $\bGhd_3$ to take the form
		\begin{equation}\label{Lyapaut}
			\left[\hy_1^2+\hy_2^2\right]\mBhl\byh,\qquad\text{where}\quad
			\mBhl\equiv\left(\begin{array}{rr} \bRl&-\bIl\\ \bIl&\bRl \end{array}\right)
		\end{equation}
		and we again call the elements of this matrix Lyapunov coefficients.
\end{itemize}
(I.e. after transformation, only a multiple of the resonant cubic terms \eqref{cubres} remains.)

After this simplification, and under Assumption~\ref{autassum1}, if we now insert $\zeta=0$ together with
\begin{subequations}\label{csubsolaut}
\begin{gather}
	\byh\tp=\bahzp,\quad\wty\tp=\zz,\quad\yo\tp=\zz,\label{csubsolauta}\\
	\lam=\lamz-\ep^2\frac{\bRlz}{\daRlz}\quad\text{and}\quad\om=\omz+\ep^2\frac{\daRlz\bIlz-\daIlz\bRlz}{\daRlz}\label{csubsolautb}
\end{gather}
\end{subequations}
into the left-hand side of \eqref{csubdeaut}, we can easily see that the remainder is $O(\ep^3)$ for \eqref{csubdeauta}, $O(\ep^2)$ for \eqref{csubdeautb}
and \eqref{csubdeautc}, and zero for \eqref{csubdeautd} and \eqref{csubdeaute}.
Consequently, by transforming \eqref{csubsolauta} back through \eqref{hhtodefaut}, i.e.
\begin{equation}\label{csubsolauttrans}
\begin{aligned}
	\bzh\tp&=\bahzp+\rep\bhh(\ep\bahzp;\theta,\lam)\\
	\wtz\tp&=\rep\wbh(\ep\bahzp;\theta,\lam)\\
	\zo\tp&=\rep\ho(\ep\bahzp;\theta,\lam),
\end{aligned}
\end{equation}
we obtain an asymptotic solution for \eqref{cmnfdeaut}. Since Theorem~\ref{TheoE3aut} already displays such a solution, i.e. $\lam^E(\ep)$, $\om^E(\ep)$,
$\eta^E(\ep)$ and
\[\bz^E\tpe\equiv\begin{pmatrix}\bzh^E\tpe\\ \wtz^E\tpe\\ \zo^E\tpe\end{pmatrix},\]
this must match with \eqref{csubsolautb} and \eqref{csubsolauttrans}.
Thus we obtain
\begin{equation}\label{cmnfsolauta}
\begin{gathered}
	\lam^E(\ep)=\lamz-\ep^2\frac{\bRlz}{\daRlz},\quad\eta^E(\ep)=\ep\kaplz,\\
	\om^E(\ep)=\omz+\ep^2\frac{\daRlz\bIlz-\daIlz\bRlz}{\daRlz}
\end{gathered}
\end{equation}
and, through \eqref{csubsolauttrans},
\begin{subequations}\label{cmnfsolautb}
\begin{align}
	\bzh^E\tpe&=\bahzp+\ep\left\{\bah_0\tlz+\bah_2\tlz\cos{2\phi}+\bbh_2\tlz\sin{2\phi}\right\}\notag\\
		&\qquad+\ep^2\left\{\bah_1\tlz\cos{\phi}+\bbh_1\tlz\sin{\phi}\right.\notag\\
		&\qquad\qquad\left.+\bah_3\tlz\cos{3\phi}+\bbh_3\tlz\sin{3\phi}\right\}\\
	\wtz^E\tpe&=\ep\left\{\wba_0\tlz+\wba_2\tlz\cos{2\phi}+\wbb_2\tlz\sin{2\phi}\right\}+O(\ep^2)\\
	\zo^E\tpe&=\ep\left\{\oa_0\tlz+\oa_2\tlz\cos{2\phi}+\ob_2\tlz\sin{2\phi}\right\}+O(\ep^2).
\end{align}
\end{subequations}
Finally, by comparing \eqref{cmnfsolautb} and \eqref{E3solaut}, we see that the coefficients of $\bhh(.;\theta,\lamz)$, $\wbh(.;\theta,\lamz)$ and
$\ho(.;\theta,\lamz)$ are given \ul{exactly} by the coefficients of the Fourier $\phi$-modes in the $\bzh_2^E\tp$ and $\bzh_3^E\tp$ terms of $\bzh^E\tpe$,
the $\wtz_2^E\tp$ term in $\wtz^E\tpe$ and the $\zo_2^E\tp$ term in $\zo^E\tpe$ for \eqref{E3solaut}. Moreover, $\kaplz=\eta_2^E$ and,
by comparing \eqref{cmnfsolauta} and \eqref{E3solaut}, the Lyapunov coefficients $\bRlz$ and $\bIlz$ in \eqref{Lyapaut} are given \ul{exactly} by
\begin{equation}\label{exactLaut}
	\bRlz=-\daRlz\lam_2^E\quad\text{and}\quad\bIlz=\om_2^E-\daIlz\lam_2^E.
\end{equation}
Thus we see that Assumption~\ref{autassum2} is equivalent to $\bRlz\neq0$.

To calculate the expansion in \eqref{E3solaut}, however, requires (through $\bG$) explicit knowledge of the second and third derivatives of $\bF$, so it
is practically much more convenient to approximate not only the coefficients of $\bhh(.;\theta,\lamz)$, $\wbh(.;\theta,\lamz)$ and $\ho(.;\theta,\lamz)$
but also the Lyapunov coefficients $\bRlz$ and $\bIlz$ by using instead the \ul{$M=3$} Fourier $\phi$-approximation in \eqref{F3solaut},
i.e. $\lam^F(\ep)$, $\om^F(\ep)$, $\eta^F(\ep)$ and
\[\bz^F\tpe\equiv\left\{
\begin{aligned}
	\bzh^F\tpe&\equiv\bah_0^F(\theta,\ep)+\sum_{m=1}^3\bah_m^F(\theta,\ep)\cos{m\phi}+\bbh_m^F(\theta,\ep)\sin{m\phi}\\
	\wtz^F\tpe&\equiv\wba_0^F(\theta,\ep)+\sum_{m=1}^3\wba_m^F(\theta,\ep)\cos{m\phi}+\wbb_m^F(\theta,\ep)\sin{m\phi}\\
	\zo^F\tpe&\equiv\oa_0^F\tep+\sum_{m=1}^3\oa_m^F\tep\cos{m\phi}+\ob_m^F\tep\sin{m\phi}.
\end{aligned}
\right.\]
\begin{theorem}\label{Theoerraut}
	Using the asymptotic error results in Theorem~\ref{TheoF3aut} on page~\pageref{TheoF3aut}, our practical approximate formulae are
	\begin{align*}
		\bRlz&=-\daRlz\frac{\lam^F(\ep)-\lamz}{\ep^2}+O(\ep^2)\\
		\bIlz&=\frac{\om^F(\ep)-\omz}{\ep^2}-\daIlz\frac{\lam^F(\ep)-\lamz}{\ep^2}+O(\ep^2)\\
		\kaplz&=\frac{\eta^F(\ep)}{\ep}+O(\ep^2),
	\end{align*} 
	together with
	{\renewcommand{\arraystretch}{1.5}
	\[\begin{array}{|c|c|}\hline
		m=0,2&m=2\\ \hline\hline
		\bah_m\tlz=\rep\bah_m^F(\theta,\ep)+O(\ep^2)&\bbh_m\tlz=\rep\bbh_m^F(\theta,\ep)+O(\ep^2)\\ \hline
		\wba_m\tlz=\rep\wba_m^F(\theta,\ep)+O(\ep^2)&\wbb_m\tlz=\rep\wbb_m^F(\theta,\ep)+O(\ep^2)\\ \hline
		\oa_m\tlz=\rep\oa_m^F(\theta,\ep)+O(\ep^2)&\ob_m\tlz=\rep\ob_m^F(\theta,\ep)+O(\ep^2)\\ \hline
	\end{array}\]}
and
	{\renewcommand{\arraystretch}{1.5}
	\[\begin{array}{|l||l|l|}\hline
		m=1,3&\bah_m\tlz=\rept\bah_m^F(\theta,\ep)+O(\ep^2)&\bbh_m\tlz=\rept\bbh_m^F(\theta,\ep)+O(\ep^2)\\ \hline
	\end{array}\]}
\end{theorem}

We conclude by remarking that the final comment in \S\ref{approxnfpf} applies here as well.

\subsection{Higher-order Fourier approximation of tori}\label{hoFaut}
In order to compute higher-order approximations for our invariant tori, we must employ a more suitable parametrisation than \eqref{torusaut1} and 
the presence of the zero Floquet exponent in $\mEz$ means that this parametrisation is different from \eqref{toruspf2}.
Thus we use the normal bundle of the approximate torus
\begin{equation}\label{torusaut2}
	\bvztl+\ep\mPtl\bazp
\end{equation}
and, in \eqref{csubdeaut},
\begin{itemize}
	\item replace $\byhtp$ with $\left[1+\rhotp\right]\bahzp$ for unknown $\rho:\So\times\So\mapsto\Rn$,
	\item replace $\yotp$ by $0$,
	\item allow both $\om:\So\times\So\mapsto\R$ and $\zeta:\So\times\So\mapsto\R$ to be unknown \ul{functions}.
\end{itemize}
This links up with the invariance condition used in \cite{Moore96} for continuation of tori, and corresponds to using polar co-ordinates in the critical
$2$-dimensional subspace. In \eqref{csubdeautd} and \eqref{csubdeaute}, there is now no need for scalar phase conditions, and the scalar amplitude
equation simplifies to a zero-mean condition for $\rho$, i.e.
\begin{equation}\label{zeromaut}
	\left\langle\left\langle\rhotp,1\right\rangle\right\rangle=0.
\end{equation}
Thus our equations for $\rho$ and $\om$ in \eqref{csubdeauta} decouple to become
\begin{subequations}\label{centreaut}
\begin{gather}
	\begin{split}
		\bahzp\cdot\bGhd\Big(&\left[1+\rhotp\right]\bahzp,\wtytp,0,\lam,\zetp,\theta;\ep\Big)\\
		&-\left[(1+\ep^2\kapl+\ep\zetp)\pdt+\omtp\pdp\right]\rhotp=0
	\end{split}\label{centreauta}\\
\intertext{and}
	\begin{split}
		\frac1{1+\rhotp}\badhzp\cdot\bGhd\Big(\left[1+\rhotp\right]\bahzp,\wtytp,0,\lam,&\zetp,\theta;\ep\Big)\\
		&-\omtp=0,
	\end{split}\label{centreautb}
\end{gather}
\end{subequations}
while the hyperbolic equations in \eqref{csubdeautb} remain
\begin{equation}\label{hyperaut}
	\begin{split}
		\wbGd\Big(\left[1+\rhotp\right]&\bahzp,\wtytp,0,\lam,\zetp,\theta;\ep\Big)\\
		&-\left[\big(1+\ep^2\kapl+\ep\zetp\big)\pdt+\omtp\pdp\right]\wtytp=\zz
	\end{split}
\end{equation}
and \eqref{csubdeautc} becomes
\begin{equation}\label{zeroaut}
	\God\Big(\left[1+\rhotp\right]\bahzp,\wtytp,0,\lam,\zetp,\theta;\ep\Big)-\left[\zetp+\ep\kapl\right]=0.
\end{equation}
The crucial leading terms in \eqref{centreaut} are
\begin{subequations}\label{ltermaut}
\begin{equation}
\begin{split}
	\bahzp\cdot\bGhd&\Big(\left[1+\rhotp\right]\bahzp,\wtytp,0,\lam,\zetp,\theta;\ep\Big)\\
	&=\left[1+\rhotp\right]\aRl+\ep^2\bRl\left[1+\rhotp\right]^3+O(\ep^3)
\end{split}
\end{equation}
and
\begin{equation}\label{ltermautb}
\begin{split}
	\frac1{1+\rhotp}\badhzp\cdot\bGhd&\Big(\left[1+\rhotp\right]\bahzp,\wtytp,0,\lam,\zetp,\theta;\ep\Big)\\
	&=\aIl+\ep^2\bIl\left[1+\rhotp\right]^2+O(\ep^3),
\end{split}
\end{equation}
while in \eqref{zeroaut} we have
\begin{equation}\label{ltermautc}
\begin{split}
	\God\Big(\left[1+\rhotp\right]\bahzp,\wtytp,&0,\lam,\zetp,\theta;\ep\Big)\\
	&=\ep\kapl\left[1+\rhotp\right]^2+O(\ep^2).
\end{split}
\end{equation}
\end{subequations}

Although we needed to introduce $\zeta$ through \eqref{defkappa} in order to obtain the correct normal form in \S\ref{approxnfaut}, it is now simpler to
describe our final system of equations in terms of
\[\eta\tp\equiv\zeta\tp+\ep\kapl.\]
Thus we can re-write \eqref{zeroaut} as
\begin{equation}\label{defeta}
	\Godd\Big(\rhotp,\wtytp,\lam,\etp,\theta;\ep\Big)-\etp=0
\end{equation}
and use \eqref{defeta} to \ul{define} $\etp$ in terms of $\lam$, $\rhotp$ and $\wtytp$ for $\lvert\ep\rvert$ sufficiently small. (Since $\God$ depends
linearly on $\zeta$ in \eqref{zeroaut}, and thus $\Godd$ depends linearly on $\eta$ in \eqref{defeta}, this is particularly simple.) Similarly, we can
re-write \eqref{centreautb} as
\begin{equation}\label{defom}
	\frac1{1+\rhotp}\badhzp\cdot\bGhdd\Big(\rhotp,\wtytp,\lam,\theta;\ep\Big)-\omtp=0
\end{equation}
by inserting $\etp$ from \eqref{defeta} into $\bGhd$; hence \eqref{defom} defines $\omtp$ in terms of $\lam$, $\rhotp$ and $\wtytp$. Finally, we can
re-write \eqref{centreauta} and \eqref{hyperaut} as
\begin{subequations}\label{fineqaut}
\begin{gather}
	\begin{split}
		\bahzp\cdot\bGhdd\Big(\rhotp,&\wtytp,\lam,\theta;\ep\Big)\\
		&-\left[\big(1+\ep\etp\big)\pdt+\omtp\pdp\right]\rhotp=0
	\end{split}\\
	\begin{split}
		\wbGdd\Big(\rhotp,&\wtytp,\lam,\theta;\ep\Big)\\
		&-\left[\big(1+\ep\etp\big)\pdt+\omtp\pdp\right]\wtytp=\zz
	\end{split}
\end{gather}
\end{subequations}
by inserting $\etp$ from \eqref{defeta} into $\bGhd$ and $\wbGd$ respectively. In \cite{Sack64,Sack65} it is proved that the system of equations
\eqref{fineqaut} and \eqref{zeromaut} has a locally unique solution $(\lam,\rho,\wty)$,
for $\lvert\ep\rvert$ sufficiently small, by considering the iteration
\begin{subequations}\label{iteraut}
\begin{gather}
	\begin{split}
		\left\{2\ep^2\bRlz-\left[\big(1+\ep\eta^{(k)}\tp\big)\pdt+\om^{(k)}\tp\pdp\right]\right\}&\rho^{(k+1)}\tp\\
		+\daRz&\delamko=\hrktp
	\end{split}\label{iterauta}\\
	\begin{split}
		\left\{\wEz-\left[\big(1+\ep\eta^{(k)}\tp\big)\pdt+\om^{(k)}\tp\pdp\right]\right\}&\wty^{(k+1)}\tp\\
		&\quad=\wbrktp
	\end{split}\\
	\left\langle\left\langle\rho^{(k+1)}\tp,1\right\rangle\right\rangle=0,
\end{gather}
\end{subequations}
where $\delamko\equiv\lamko-\lamk$ and
\begin{align*}
	\hrktp&\equiv 2\ep^2\bRlz\rho^{(k)}\tp-\bahzp\cdot\bGhdd\Big(\rho^{(k)}\tp,\wty^{(k)}\tp,\lamk,\theta;\ep\Big)\\
	\wbrktp&\equiv\wEz\wty^{(k)}\tp-\wbGdd\Big(\rho^{(k)}\tp,\wty^{(k)}\tp,\lamk,\theta;\ep\Big),
\end{align*}
with starting values
\[\rho^{(0)}\tp=0,\;\wty^{(0)}\tp=\zz,\;\lams=\lamz.\]
(Note that $\eta^{(k)}\tp$ and $\om^{(k)}\tp$ are defined through \eqref{defeta} and \eqref{defom} respectively, using the values
$\lamk$, $\rho^{(k)}\tp$ and $\wty^{(k)}\tp$.)
The key idea behind showing that these iterates remain bounded
and then converge is to integrate \eqref{iterauta} against $\rho^{(k+1)}\tp$, after which the left-hand side becomes
\[\left\langle\left\langle\left\{2\ep^2\bRlz+\ep\pdt\eta^{(k)}\tp+\pdp\om^{(k)}\tp\right\}\rho^{(k+1)}\tp,\rho^{(k+1)}\tp\right\rangle\right\rangle.\]
Since \eqref{ltermautb} shows that the leading \ul{non-constant} term in $\om^{(k)}\tp$ is $O(\ep^3)$,
and \eqref{ltermautc} together with \eqref{zeroaut} shows that the leading \ul{non-constant} term in $\eta^{(k)}\tp$ is $O(\ep^2)$,
Assumption~\ref{autassum2} ensures that the last expression is a definite quadratic term in $\rho^{(k+1)}\tp$ for $\lvert\ep\rvert$ sufficiently small 
and this is sufficient for \cite{Sack64} to prove the following theorem.
\begin{theorem}\label{Sackaut}
	Suppose $\bF$ in \eqref{autde} has $r\ge5$ continuous derivatives for $(\lam,\bx)$ in a neighbourhood of $(\lamz,\bvzt)$. Then $\exists\,{\ep_r>0}$
	such that for $\lvert\ep\rvert<\ep_r$ \eqref{fineqaut} has a locally unique solution
	\[\lamz(\ep),\;\rhoz\tpe,\;\wtyz\tpe\]
	with $\rhoz(.,.;\ep),\;\wtyz(.,.;\ep)$ having $(r-1)$ Lipchitz continuous derivatives. This means that both $\byh$ and $\bzh,\;\wtz$
	(through \eqref{hhtodefaut}) have this degree of smoothness, and so, through \eqref{torusaut1}, do the invariant tori as manifolds.
\end{theorem}

\noindent As in Theorem~\ref{Sackpf}, in general $\ep_r\to 0$ as $r\to\infty$ and we cannot expect analytic tori.

In practice we seek an approximate solution of \eqref{fineqaut} in the form
\[\begin{aligned}
	\rhoLM\tp&\equiv\sum_{\ell,m}\rho_{\ell,m}\ee^{\ii(\ell\theta+m\phi)}\\
	\wtyLM\tp&\equiv\sum_{\ell,m}\wty_{\ell,m}\ee^{\ii(\ell\theta+m\phi)}
\end{aligned}\qquad\lvert\ell\rvert\le L,\;\lvert m\rvert\le M\]
with $\rho_{-\ell,-m},\;\wty_{-\ell,-m}$ the conjugates of $\rho_{\ell,m},\;\wty_{\ell,m}$. These functions must satisfy
\begin{subequations}\label{fineqFaut}
\begin{gather}
\begin{split}
	&\left\langle\left\langle\bahzp\cdot\bGhdd\Big(\right.\right.\rhoLM\tp,\wtyLM\tp,\lam,\theta;\ep\Big)\\
	&\quad\left.\left.-\left[\big(1+\ep\etp\big)\pdt+\omtp\pdp\right]\rhoLM\tp,\ee^{\ii(\ell\theta+m\phi)}\right\rangle\right\rangle=0
\end{split}\\
\begin{split}
	\left\langle\left\langle\wbGdd\Big(\right.\right.&\rhoLM\tp,\wtyLM\tp,\lam,\theta;\ep\Big)\\
	&\left.\left.-\left[\big(1+\ep\etp\big)\pdt+\omtp\pdp\right]\wtyLM\tp,\bx\ee^{\ii(\ell\theta+m\phi)}\right\rangle\right\rangle=\zz
\end{split}\\
	\rho_{0,0}=0,
\end{gather}
\end{subequations}
for $\lvert\ell\rvert\le L$, $\lvert m\rvert\le M$ and $\forall{\bx\in\Rnth}$, with $\eta^{(k)}\tp$ and $\om^{(k)}\tp$ defined 
through \eqref{defeta} and \eqref{defom} respectively, using the values $\lam$, $\rhoLM\tp$ and $\wtyLM\tp$.
As shown in \cite{Sam91}, the analogous iteration to \eqref{iteraut} also converges here and gives the following theorem.
\begin{theorem}\label{Samaut}
	Under the conditions of Theorem~\ref{Sackaut}, \eqref{fineqFaut} has a locally unique solution
	\[\lamLM^F(\ep),\;\rhoLM^F\tpe,\;\wtyLM^F\tpe\]
	satisfying
	\[\max\left\{\begin{array}{c}
		\lvert\lamLM^F(\ep)-\lamz(\ep)\rvert\\
		\lVert\rhoLM^F(.,.;\ep)-\mSLM\rhoz(.,.;\ep)\rVert_{L_2}\\
		\lVert\wtyLM^F(.,.;\ep)-\mSLM\wtyz(.,.;\ep)\rVert_{L_2}
	\end{array}\right\}\le C
	\max\left\{\begin{array}{c}
		\lVert\left(\mI-\mSLM\right)\rhoz(.,.;\ep)\rVert_{H^1}\\
		\lVert\left(\mI-\mSLM\right)\wtyz(.,.;\ep)\rVert_{H^1}
	\end{array}\right\}.\]
\end{theorem}
As in \S\ref{hoFpf}, we comment on the implementation of this algorithm in \S\ref{Conclu}.

\subsection{Numerical results}\label{numresaut}
We consider a numerical example for which a group orbit structure leads to an interesting simplification of the general Neimark--Sacker bifurcation
equations: this is the Kuramoto--Sivashinsky equation in the form
\begin{equation}\label{ksequ}
	\pudt\xt=-4\pfudx\xt-\lam\left[\ptudx\xt+\thalf\pdx\left(u^2\xt\right)\right],
\end{equation}
with $u$ being both $\tpi$-periodic and having zero mean in $x$ \cite{KNS90,Tem88}. We immediately obtain a finite-dimensional autonomous system by
restricting to the Fourier approximation
\begin{equation}\label{ksmodes}
	u\xt\approx\sum_{\ell=-L}^L u_{\ell}(t)\ee^{\ii\ell x}\qquad\begin{aligned}u_0(t)&=0\\ u_{-\ell}(t)&=\ol{u}_{\ell}(t)\end{aligned},
\end{equation}
and making use of the conjugacy condition leads to the complex system
\begin{equation}\label{ksL}
	\dbudt=-4\mDL^4\bu+\lam\left[\mDL^2\bu-\ii\mDL\bQC(\bu)\right]\qquad\bu\in\CL,	
\end{equation}
where $\mDL$ is the $L\times L$ diagonal matrix with entries $1,\dots,L$ and the quadratic function $\bQC:\CL\mapsto\CL$ is defined by
\begin{equation}\label{kserr}
	\left[\QC(\bu)\right]_{\ell}\equiv\thalf\sum_{j=1}^{\ell-1}u_{\ell-j}u_j+\sum_{j=1}^{L-\ell}\ol{u}_j u_{\ell+j}.
\end{equation}
(Thus \eqref{kserr} is our only discretisation error.)
We describe below the sequence of computations which leads to Neimark--Sacker bifurcation for \eqref{ksL}:
this mimics some of the numerical results in \cite{KNS90}, which should be referred to for further information.
These computations exhibit our fundamental Crandall--Rabinowitz formulation in three different bifurcation situations.

\paragraph{a) Bifurcation from the trivial solution} \eqref{ksL} has the trivial stationary solution curve $\bu\equiv 0\;\forall\lam$,
which is stable for $\lam<4$, and nontrivial stationary solutions bifurcate at
\begin{equation}\label{kstbif}
	\lamztb=4\ell^2\quad\ell=1,2,\dots,L.
\end{equation}
These nontrivial stationary solutions are \ul{not} isolated, since the autonomous nature of \eqref{ksequ} implies that if $u\xt$ is a solution then
so is $u(x+\alpha,t)\;\forall\alpha\in\R$. Consequently, in order to apply the \IFT\ and Newton's method, we must eliminate this multiplicity by
either imposing a phase condition or a symmetry restriction. Since we are interested in the first bifurcation branch, i.e. $\lamztb=4$ in
\eqref{kstbif}, it is simplest to consider only stationary solutions of the form
\[\bu=\ii\bs\quad\text{with}\quad\bs\in\RL\]
for \eqref{ksL}: this leads to the real system
\begin{equation}\label{kssymm}
	-4\mDL^4\bs+\lam\left[\mDL^2\bs+\mDL\bQI(\bs)\right]=\zz,
\end{equation}
where the quadratic function $\bQI:\RL\mapsto\RL$ is defined by
\[\left[\QI(\bs)\right]_{\ell}\equiv\thalf\sum_{j=1}^{\ell-1}s_{\ell-j}s_j-\sum_{j=1}^{L-\ell}s_j s_{\ell+j}.\]
For small $\lvert\ep\rvert$, we move onto the bifurcating curve of nontrivial stationary solutions by seeking solutions of \eqref{kssymm} in the
Crandall--Rabinowitz formulation $\bs\equiv\ep\hs$, with amplitude condition $\hat{s}_1=1$. Hence, with starting values
\[\lams=4,\quad\hsf=\bef,\]
the iteration in \cite{CR71} can be written
\[\lamko=4-\ep r_1^{(k)},\qquad \hat{s}_{\ell}^{(k+1)}=\ep r_{\ell}^{(k)}/(4\ell^3-\ell\lamko)\quad\ell=2,\dots,L,\]
where
\[\brk\equiv\lamk\bQI(\hsk).\]

\paragraph{b) Continuation of stationary solutions} Having moved away from the bifurcation point at $\lamztb=4$, we can follow
the branch of nontrivial stationary solutions by applying a standard continuation algorithm \cite{AG90} to \eqref{kssymm}. This branch is always
parametrisable by $\lam$, and so we can refer to solutions of \eqref{kssymm} by $(\lam,\bs(\lam))$ and the Jacobian matrix at solutions by
\begin{equation}\label{kssymmp}
	\mJIss(\lam)\equiv-4\mDL^4+\lam\left[\mDL^2+\mDL\mTa(\bs(\lam))-\mDL\mH(\bs(\lam))\right],
\end{equation}
where in Matlab notation
\begin{gather*}
	\mTa(\bs)\equiv\text{\texttt{toeplitz}}([0\; s(1:L-1)], [0\; -s(1:L-1)])\\
	\mH(\bs)\equiv\text{\texttt{hankel}}([s(2:L)\; 0]).
\end{gather*}
The eigenvalues of \eqref{kssymmp} remain strictly in the left-half plane but this matrix, however, only measures the effect of symmetric perturbations.
To consider the effect of symmetry-breaking perturbations we must monitor the matrix
\begin{equation}\label{kssymmb}
	\mJRss(\lam)\equiv-4\mDL^4+\lam\left[\mDL^2+\mDL\mTa(\bs(\lam))+\mDL\mH(\bs(\lam))\right],
\end{equation}
which always has a null-vector $\mDL\bs(\lam)$ because symmetry was imposed specifically to eliminate non-isolated stationary solutions. As $\lam$ moves
away from $\lamztb=4$, all the other $L-1$ eigenvalues of \eqref{kssymmb} remain at first strictly in the left-half plane but, as $\lam$ approaches
$\lamzrw\approx 13$, our zero eigenvalue becomes defective, with algebraic multiplicity two. At this value of $\lam$, we denote the null-vector of
\eqref{kssymmb} by $\berw$, with normalisation $\lVert\berw\rVert=1$, and the generalised eigenvector by $\bsigrw$, with normalisation
$(\czrw)^2+\lVert\bsigrw\rVert^2=1$, where
\[\mJRss(\lamzrw)\bsigrw=\czrw\mDL\bs(\lamzrw)\qquad\text{with $\left(\berw\right)^T\bsigrw=0$.}\]
As part of our continuation algorithm, we can monitor the real part of the eigenvalues of \eqref{kssymmb} and detect a crossing of the imaginary
axis: a simple secant iteration then accurately determines the value of $\lam$ at which bifurcation occurs and this is displayed in Figure~\ref{ksbifpts}.
\begin{figure}[t]\begin{center}
	\begin{minipage}{\textwidth}
		\[\begin{array}{|c||c|c|c|c|}\hline
			L&\lamzrw&\czrw&\lamzns&\omzns\\ \hline
			8  &13.0038442196&  0.9990409957&  17.3973078781&  3.3475479311\\
			16  &13.0038442196&  0.9990409957&  17.3973072209&  3.3475479124\\ \hline
		\end{array}\]
	\end{minipage}\end{center}
	\caption{\label{ksbifpts}Approximation of bifurcation points for the Kuramoto--Sivashinsky equation}
\end{figure}
We can also check that the crossing is transversal, by using a simple 2nd-order centered finite difference (with step $h$) to obtain the
following approximations to the critical eigenvalue derivative.
\[\begin{array}{|c||c|c|c|}\hline
	h&0.1&0.01&0.001\\ \hline
	\text{Eigenvalue speed}  &6.127497&  6.127414&  6.127418\\ \hline
\end{array}\]

\paragraph{c) Bifurcation to rotating waves} This loss of stability is associated with the creation of a special type of periodic orbit called a
\ul{rotating wave}. It is a solution of \eqref{ksequ} with \eqref{ksmodes} having the form
\[u_{\ell}(t)\equiv u_{\ell}\:\ee^{\ii\ell ct},\]
where the unknown wave-speed $c\in\R$ plays the role of ``frequency''. The important practical point is that these rotating waves are as easy to compute as
stationary solutions, since under the moving frame
\[\xi\equiv x+ct\]
they satisfy
\begin{equation}\label{ksrw}
	-4\dfvdxi-\lam\left[\dtvdxi+\thalf\ddxi\left(v^2\right)\right]-c\dvdxi=0,	
\end{equation}
where now $v$ is $\tpi$-periodic and has zero mean in $\xi$. Hence, instead of \eqref{ksmodes}, we use
\begin{equation}\label{ksrwmod}
	v(\xi)\approx\sum_{\ell=-L}^L v_{\ell}\ee^{\ii\ell\xi}\qquad\begin{aligned}v_0&=0\\ v_{-\ell}&=\ol{v}_{\ell}\end{aligned}
\end{equation}
and arrive at the complex system
\begin{equation}\label{ksrwL}
	-4\mDL^4\bv+\lam\left[\mDL^2\bv-\ii\mDL\bQC(\bv)\right]-\ii c\mDL\bv=\zz\qquad\bv\in\CL,	
\end{equation}
which is the analogue of \eqref{ksL}.
We can then move onto the curve of rotating waves by seeking a solution of \eqref{ksrwL} in the Crandall--Rabinowitz formulation
\[c\equiv\ep\hc\quad\text{and}\quad\bv\equiv\ii\bs(\lam)+\ep\left[\bvR+\ii\bvI\right]\qquad\bvR,\bvI\in\RL\]
for small $\lvert\ep\rvert$. Just as for ordinary Hopf bifurcation, we must complement \eqref{ksrwL} with amplitude and phase conditions, and these are
\[\left(\bsigrw\right)^T\bvR+\czrw\hc=1\quad\text{and}\quad\left(\berw\right)^T\bvR=0.\]
Thus, splitting \eqref{ksrwL} into real and imaginary parts, our analogue of the Hopf bifurcation iteration in section \S\ref{Hopfbif} is
\[\begin{bmatrix}
	\mJRss(\lamk)&\bd_1^{(k)}&\bd_2^{(k)}\\
	\left(\berw\right)^T&0&0\\
	\left(\bsigrw\right)^T&\czrw&0
	\end{bmatrix}
	\begin{bmatrix}\bvRko\\ \hcko\\ \delamko\end{bmatrix}=\begin{bmatrix}\ep\brRk\\0\\ 1\end{bmatrix}
	\quad\text{and}\quad\mJIss(\lamk)\bvIko=\ep\brIk,\]
where
\begin{align*}
	\bd_1^{(k)}&\equiv\mDL\bs(\lamk)\quad\text{and}\quad\bd_2^{(k)}\equiv(\mJRss)'(\lamzrw)\bvRk+\hck\mDL\bs'(\lamzrw)\\
	\brRk&\equiv-\hck\mDL\bvIk-\lamk\mDL\Imag\left\{\bQC(\bvRk+\ii\bvIk)\right\}\\
	\brIk&\equiv\mathrel{\phantom{-}}\hck\mDL\bvRk+\lamk\mDL\Real\left\{\bQC(\bvRk+\ii\bvIk)\right\},
\end{align*}
with starting values
\[\lams=\lamzrw,\;\hcs=\czrw,\;\bvRs=\bsigrw,\;\bvIs=\zz.\]

\begin{figure}[t]
	\begin{center}
		\includegraphics[width=.9\textwidth]{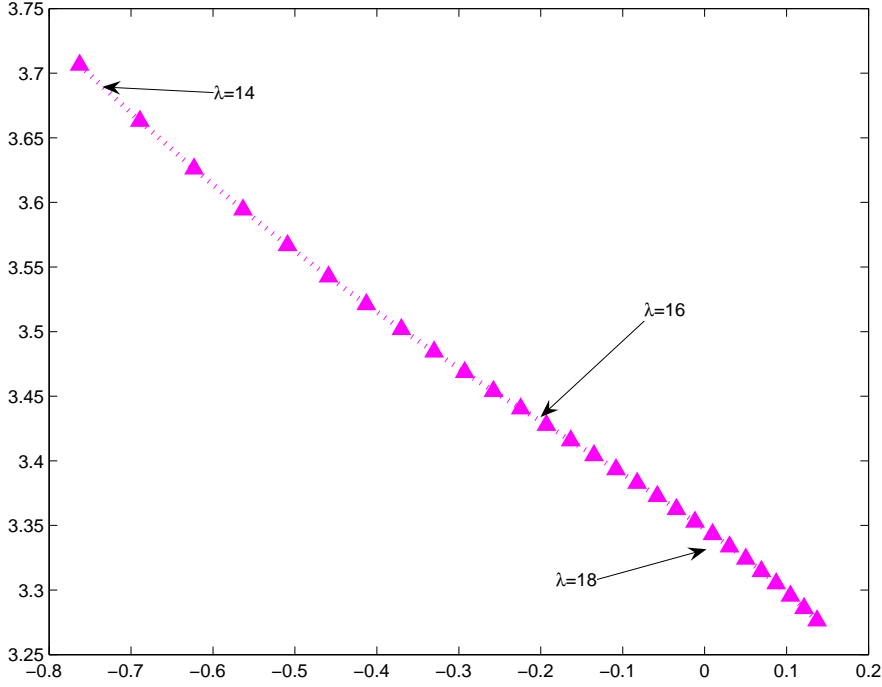}
	\end{center}
	\caption{Movement in $\C$ for critical eigenvalue of rotating waves\label{ksevalfig}}
\end{figure}

\paragraph{d) Continuation of rotating waves} Having moved away from this pseudo-Hopf bifurcation point $\lamzrw$, we can follow the branch of rotating
waves by applying a standard continuation algorithm \cite{AG90} to \eqref{ksrwL}. This branch is parametrisable by $\lam$ and so, splitting $\bv$ into
real and imaginary parts, we can refer to the solutions of \eqref{ksrwL} by $(\lam, c(\lam), \bvR(\lam)+\ii\bvI(\lam))$.
The analogue of Floquet exponents for the rotating waves are the eigenvalues of
\begin{equation}\label{ksFloqu}
\mJrw(\lam)\equiv\begin{bmatrix}\mJRRrw(\lam)&\mJRIrw(\lam)\\ \mJIRrw(\lam)&\mJIIrw(\lam)\end{bmatrix}\in\R^{{2L}\times{2L}},
\end{equation}
where
\begin{align*}
	\mJRRrw(\lam)&\equiv\frac1{c(\lam)}\left\{-4\mDL^4+\lam\mDL^2+\lam\mDL\left[\mTa(\bvI(\lam))+\mH(\bvI(\lam))\right]\right\}\\
	\mJRIrw(\lam)&\equiv\frac{\lam}{c(\lam)}\mDL\left[\mTs(\bvR(\lam))-\mH(\bvR(\lam))\right]+\mDL\\
	\mJIRrw(\lam)&\equiv\frac{-\lam}{c(\lam)}\mDL\left[\mTs(\bvR(\lam))+\mH(\bvR(\lam))\right]-\mDL\\
	\mJIIrw(\lam)&\equiv\frac1{c(\lam)}\left\{-4\mDL^4+\lam\mDL^2+\lam\mDL\left[\mTa(\bvI(\lam))-\mH(\bvI(\lam))\right]\right\}
\end{align*}
and in Matlab notation
\[\mTs(\bvR)\equiv\text{\texttt{toeplitz}}([0\; v_{\scriptscriptstyle\mathrm{R}}(1:L-1)]).\]
As with periodic orbits, one of these is always zero since
\[\berw(\lam)=\begin{bmatrix}\beRrw(\lam)\\ \beIrw(\lam)\end{bmatrix}\equiv\begin{bmatrix}-\mDL\bvI(\lam)\\ \hphantom{-}\mDL\bvR(\lam)\end{bmatrix}\]
is a null-vector because of the autonomous nature of \eqref{ksrw}. Thus \eqref{ksrwL} must be complemented by a phase condition
\[\left(\beRrw(\lam_\mathrm{prev})\right)^T\bvR+\left(\beIrw(\lam_\mathrm{prev})\right)^T\bvI=0,\] 
where $\berw(\lam_\mathrm{prev})$ is obtained from the solution at the previous value of $\lam$.
Apart from this, all the other $2L-1$ eigenvalues of \eqref{ksFloqu} lie strictly in the left-half plane until $\lam$ approaches $\lamzns\approx 17.4$,
when a complex-conjugate pair $\pm\ii\omzns$ cross the imaginary axis with the complex eigenvector satisfying
\[\mJrw(\lamzns)\begin{bmatrix}\bsigRns\\ \bsigIns\end{bmatrix}=\ii\omzns\begin{bmatrix}\bsigRns\\ \bsigIns\end{bmatrix}
	\qquad\text{for $\bsigRns,\bsigIns\in\CL$}.\]
As part of our continuation algorithm, we can monitor the real part of the eigenvalues of \eqref{ksFloqu} and detect a crossing of the imaginary
axis: a simple secant iteration then accurately determines the value of $\lam$ at which bifurcation occurs, with
\[\bens=\begin{bmatrix}\beRns\\ \beIns\end{bmatrix}\equiv \begin{bmatrix}\beRrw(\lamzns)\\ \beIrw(\lamzns)\end{bmatrix}\]
denoting the null-vector there. The variation with $\lam$ (in the upper-half of the complex plane) of this critical complex-conjugate eigenvalue is shown
in Figure~\ref{ksevalfig}, which may be compared with Figure 4.2 in \cite{KNS90},
while numerical values for Neimark--Sacker bifurcation are displayed in Figure~\ref{ksbifpts}. (For $L=32$, all results agreed to $10$ decimal places.)
As with bifurcation to rotating waves, we can also check that the crossing is transversal by calculating the
following approximations to the real part of the critical eigenvalue derivative.
\[\begin{array}{|c||c|c|c|}\hline
	h&0.1&0.01&0.001\\ \hline
	\text{Eigenvalue speed}  &0.0966405&  0.0966218&  0.0966216\\ \hline
\end{array}\]

\paragraph{e) Bifurcation to invariant tori} We seek invariant tori in the Crandall-Rabinowitz formulation
\[\cT(\xi,\phi)\equiv\sum_{\ell=-L}^L\left\{v_{\ell}(\lam)\ee^{\ii\ell\xi}+\ep\sum_{m=-M}^M z_{\ell,m}\ee^{\ii(\ell\xi+m\phi)}\right\}\qquad
	\begin{aligned}z_{0,m}&=0\quad\forall m\\ z_{-\ell,-m}&=\ol{z}_{\ell,m}\quad\forall{\ell,m}\end{aligned}\]
for small $\lvert\ep\rvert$, so that $\cT$, together with $\om,\eta\in\R$, solves the finite-dimensional restriction of
\[\frac1{c(\lam)}\left\{-4\pfTdxi-\lam\left[\ptTdxi+\thalf\pdxi\left(\cT^2\right)\right]\right\}-\left[1+\ep\eta\right]\pTdxi-\om\pTdp=0.\]
Here $(\lam,c(\lam),\bvR(\lam)+\ii\bvI(\lam))$ satisfy \eqref{ksrwL} and, if we denote the $\phi$-modes of $z$ by $\bz(m)\in\CL$, then by conjugacy we need
only solve for $\bz(0)\equiv\bzR(0)+\ii\bzI(0)$ with $\bzR(0),\bzI(0)\in\RL$ and $\bzp(m),\bzm(m)\in\CL$ for $m=1,\dots,M$. (Note that $\bzp$
contains the $\xi$-modes for $\ell\ge1$ and $\bzm$ for $\ell\le-1$.) Hence, dividing through by $\ep$, we may write the finite-dimensional restriction as
\[\begin{bmatrix}
	\mJRRrw(\lam)&\mJRIrw(\lam)&\hphantom{-}\mDL\bvI(\lam)\\
	\mJIRrw(\lam)&\mJIIrw(\lam)&-\mDL\bvR(\lam)
\end{bmatrix} \begin{bmatrix}\bzR(0)\\ \bzI(0)\\ \eta\end{bmatrix}=\ep\begin{bmatrix}\brR(0)\\ \brI(0)\end{bmatrix}\]
for $m=0$ and
\[\left[\mJCrw(\lam)-\ii\om m\mI\right]
	\begin{bmatrix}\bzp(m)\\ \bzm(m)\end{bmatrix}=\ep\begin{bmatrix}\brp(m)\\ \brm(m)\end{bmatrix}\]
for $m\ge1$, where
\[\mJCrw(\lam)\equiv\begin{bmatrix}\mJpprw(\lam)&\mJpmrw(\lam)\\ \mJmprw(\lam)&\mJmmrw(\lam)\end{bmatrix}\]
is the complex version of \eqref{ksFloqu} defined by
\begin{align*}
	\mJpprw(\lam)&\equiv\frac1{c(\lam)}\left\{-4\mDL^4+\lam\mDL^2-\ii\lam\mDL\mTs(\bvR(\lam)-\ii\bvI(\lam))\right\}-\ii\mDL\\
	\mJpmrw(\lam)&\equiv-\ii\frac{\lam}{c(\lam)}\mDL\mH(\bvR(\lam)+\ii\bvI(\lam))\\
	\mJmprw(\lam)&\equiv\ii\frac{\lam}{c(\lam)}\mDL\mH(\bvR(\lam)-\ii\bvI(\lam))\\
	\mJmmrw(\lam)&\equiv\frac1{c(\lam)}\left\{-4\mDL^4+\lam\mDL^2+\ii\lam\mDL\mTs(\bvR(\lam)+\ii\bvI(\lam))\right\}+\ii\mDL.
\end{align*}
The right-hand sides are defined by
\[\begin{bmatrix}\brR(0)\\ \brI(0)\end{bmatrix}\equiv
	\frac{\lam}{c(\lam)}\begin{bmatrix}-\mDL\bqI(0)\\ \hphantom{-}\mDL\bqR(0)\end{bmatrix}+
		\eta\begin{bmatrix}-\mDL\bzI(0)\\ \hphantom{-}\mDL\bzR(0)\end{bmatrix}\]
and
\[\begin{bmatrix}\brp(m)\\ \brm(m)\end{bmatrix}\equiv
	\ii\frac{\lam}{c(\lam)}\begin{bmatrix}\hphantom{-}\mDL\bqp(m)\\-\mDL\bqm(m)\end{bmatrix}+
		\ii\eta\begin{bmatrix}\hphantom{-}\mDL\bzp(m)\\-\mDL\bzm(m)\end{bmatrix},\]
with $\bqR(0),\bqI(0)\in\RL$ and $\bqp(m),\bqm(m)\in\CL$ being derived from the quadratic term
\[\thalf\mSLM\left(\sum_{\genfrac{}{}{0pt}{3}{m=-M}{\ell=-L}}^{\scriptscriptstyle L,M} z_{\ell,m}\ee^{\ii(\ell\xi+m\phi)}\right)^2
	=\sum_{\genfrac{}{}{0pt}{3}{m=-M}{\ell=-L}}^{\scriptscriptstyle L,M} q_{\ell,m}\ee^{\ii(\ell\xi+m\phi)}\qquad
	\begin{aligned}q_{0,m}&=0\quad\forall m\\ q_{-\ell,-m}&=\ol{q}_{\ell,m}\quad\forall{\ell,m}\end{aligned}\]
in the same way as their analogues for $\bz$.  We also have the phase condition
\[\left(\beRns\right)^T\bzR(0)+\left(\beIns\right)^T\bzI(0)=0\]
and, since 
\[\left[\mJCrw(\lamzns)-\ii\omz\mI\right]
	\begin{bmatrix}\bsigRns+\ii\bsigIns\\ \bsigRns-\ii\bsigIns\end{bmatrix}=\begin{bmatrix}\zz\\ \zz\end{bmatrix},\]
we can define our amplitude and phase conditions using
\[\begin{bmatrix}\bsigpns\\ \bsigmns\end{bmatrix}\equiv\alpha\begin{bmatrix}\bsigRns+\ii\bsigIns\\ \bsigRns-\ii\bsigIns\end{bmatrix},\]
where $\alpha\in\C$ is chosen so that
\[\lVert\bsigpns\rVert^2+\lVert\bsigmns\rVert^2=1.\]
Thus our Neimark--Sacker bifurcation iteration is
\[\begin{bmatrix}
	\mJRRrw(\lamk)&\mJRIrw(\lamk)&\hphantom{-}\mDL\bvI(\lamk)\\
	\mJIRrw(\lamk)&\mJIIrw(\lamk)&-\mDL\bvR(\lamk)\\
	\left[\beRns\right]^T&\left[\beIns\right]^T&0
\end{bmatrix} \begin{bmatrix}\bzRko(0)\\ \bzIko(0)\\ \etako\end{bmatrix}=\\
	=\ep\begin{bmatrix}\brRk(0)\\ \brIk(0)\\0\end{bmatrix}\]
for $m=0$;
\[\begin{bmatrix}
	\mJpprw(\lamk)-\ii\omk\mI&\mJpmrw(\lamk)&-\ii\mDL\bzpk(1)&\bd_{\scriptscriptstyle +}^{(k)}\\
	\mJmprw(\lamk)&\mJmmrw(\lamk)-\ii\omk\mI&-\ii\mDL\bzmk(1)&\bd_{\scriptscriptstyle -}^{(k)}\\
	(\bsigpns)^{\star}&(\bsigmns)^{\star}&0&0
\end{bmatrix} \begin{bmatrix}\bzpko(1)\\ \bzmko(1)\\ \delomko\\ \delamko\end{bmatrix}
	=\begin{bmatrix}\ep\brpk(1)\\ \ep\brmk(1)\\1\end{bmatrix}\]
for $m=1$; and
\[\left[\mJCrw(\lamk)-\ii m\omk\mI\right]\begin{bmatrix}\bzpko(m)\\ \bzmko(m)\end{bmatrix}=\ep\begin{bmatrix}\brpk(m)\\ \brmk(m)\end{bmatrix}\]
for $m\ge2$. Here
\[\begin{bmatrix}\brRk(0)\\ \brIk(0)\end{bmatrix}\equiv
	\frac{\lamk}{c(\lamk)}\begin{bmatrix}-\mDL\bqIk(0)\\ \hphantom{-}\mDL\bqRk(0)\end{bmatrix}+
		\etak\begin{bmatrix}-\mDL\bzIk(0)\\ \hphantom{-}\mDL\bzRk(0)\end{bmatrix},\]
\[\begin{bmatrix}\brpk(m)\\ \brmk(m)\end{bmatrix}\equiv
	\ii\frac{\lamk}{c(\lamk)}\begin{bmatrix}\hphantom{-}\mDL\bqpk(m)\\-\mDL\bqmk(m)\end{bmatrix}+
		\etak\begin{bmatrix}\hphantom{-}\mDL\bzpk(m)\\-\mDL\bzmk(m)\end{bmatrix},\]
\[\begin{bmatrix}\bd_{\scriptscriptstyle +}^{(k)}\\ \bd_{\scriptscriptstyle -}^{(k)}\end{bmatrix}\equiv
	\left[\mJCrw\right]'(\lamzns)\begin{bmatrix}\bzpk(1)\\ \bzmk(1)\end{bmatrix}.\]
and we note that, for $m=1$, our \ul{two} extra real unknowns $\om,\lam$ are compensated by \ul{one} extra complex condition.
Our starting values are
\[\lams=\lamzns,\;\oms=\omzns,\;\etas=0,\;\bzps(1)=\bsigpns,\;\bzms(1)=\bsigmns;\]
with all other components of $\bzs$ zero.
\begin{figure}[ht]\begin{center}
	\begin{minipage}{\textwidth}
		\[\begin{array}{|c|c|c|c|}\hline
			\eta&\om&\lam&c(\lam)\\\hline
			-0.1652838896&3.1325935189&17.3749646320&16.7931284737\\ \hline
		\end{array}\]
		\[\begin{array}{|c|c|}\hline
			\lVert\bz(0)\rVert&0.25\\ \hline
			\lVert\bz(1)\rVert&1.00\\ \hline
			\lVert\bz(2)\rVert&0.10\\ \hline
		\end{array}\quad
		\begin{array}{|c|c|}\hline
			\lVert\bz(3)\rVert&0.11\times 10^{-1}\\ \hline
			\lVert\bz(4)\rVert&0.13\times 10^{-2}\\ \hline
			\lVert\bz(5)\rVert&0.15\times 10^{-3}\\ \hline
		\end{array}\quad
		\begin{array}{|c|c|}\hline
			\lVert\bz(6)\rVert&0.18\times 10^{-4}\\ \hline
			\lVert\bz(7)\rVert&0.21\times 10^{-5}\\ \hline
			\lVert\bz(8)\rVert&0.25\times 10^{-6}\\ \hline
		\end{array}\]
	\end{minipage}\end{center}
	\caption{\label{kstorifig}Decay of Fourier $\phi$-modes for invariant torus at $\ep=0.5$}
\end{figure}
In Figure~\ref{kstorifig} we display the numerical results for $L=16$ and $M=8$, in particular verifying the decay of the size of the $\phi$-modes
for $\bz^T\equiv\left(\bzp^T,\bzm^T\right)$.
\paragraph{f) Concluding remarks} Finally, we emphasise the simplification in the Neimark--Sacker algorithm that the group orbit structure of the
Kuramoto--Sivashinsky equation allows. Just as the rotating waves are really periodic orbits that can be calculated as simply as stationary solutions,
so the invariant tori can be calculated as simply as periodic orbits: i.e. there is only \ul{one} explicit independent periodic variable and thus no
resonance can occur. This means that we can utilise the simple parametrisation for the tori in \S\ref{secCRformaut}
(as above with constant $\eta$ and $\om$) for an arbitrary number of $\phi$-modes, rather than being limited to $M=3$ by Assumption~\ref{autassum1}.

\section{Conclusion}\label{Conclu}
In \S\ref{Intro} we stated that \emph{the fundamental idea behind the present paper is to use the approach in \cite{Sack64} \dots
to develop a practical computational algorithm for Neimark--Sacker bifurcation}. We claim to have achieved this goal, but the final implementation of the
algorithms in \S\ref{hoFpf} and \S\ref{hoFaut} will be explored elsewhere. The two main reasons for this are the length of the present paper and the
belief that these practical questions are best-suited to a separate paper. We emphasise, however, the two key points that an
efficient algorithm must address.
\begin{enumerate}
	\item[a)] The extraction of normal form information from the simple low-order Fourier approximations in \S\ref{approxnfpf} and
		\S\ref{approxnfaut}: which then allows us to introduce the essential, but more complex, parametrisations in \S\ref{hoFpf} and
		\S\ref{hoFaut}.
	\item[b)] Our final iterations in \eqref{iterpf}, \eqref{fineqFpf}, \eqref{iteraut}  and \eqref{fineqFaut}  necessarily rely on the
		solution of linear \emph{variable-coefficient} differential equations. This raises the question of computational efficiency since,
		throughout this paper, we have utilised the mode-decoupling property for Fourier approximations of
		\emph{constant-coefficient} systems. Our solution to this problem is to make use of the precise structure of the variable-coefficient
		equations in order to \ul{pre-condition} them by suitable constant-coefficient operators \cite{Boy01,CHQZ88,Vorst03}.
\end{enumerate}
Finally, we remark on several other points which, for the sake of simplicity, were omitted earlier.
\begin{itemize}
	\item In \S\ref{NSbifpf} and \S\ref{NSbifaut} we assumed that the basic periodic orbit $\bvzt$ was known exactly. In practice, of course,
		we would have a sufficiently accurate Fourier approximation, as in \S\ref{numrespf} and \S\ref{numresaut}.
	\item In \S\ref{NSbifpf} and \S\ref{NSbifaut} we assumed that $\nm=0$ for the Floquet theory described in \S\ref{CompFT}.
		The case $\nm>0$ introduces no practical difficulties, whether these eigenvalues occur in $\mEhz$ or $\wEz$.
		In both cases, the strategy in \cite{Moore05} can be followed.
	\item We have avoided any discussion of aliasing, numerical quadrature and the \texttt{FFT} for our Fourier spectral methods \cite{Boy01,CHQZ88},
		by implicitly assuming that all integration was performed exactly. The only practical difference is that some of our errors in
		Theorems~\ref{Hopfpracerr}, \ref{Theoerrpf} and \ref{Theoerraut} may be $O(\ep)$ rather than $O(\ep^2)$. This is still sufficient
		for our purposes, but may be avoided if desired: such questions will be addressed in the future paper mentioned above.
	\item We have merely stated the smooth invariant subspace decompositions required in \eqref{Hopfinvar}, \eqref{pfinvar} and \eqref{autinvar}.
		Further information may be found in \cite{DLF00}.
\end{itemize}

\bibliographystyle{siam}
\bibliography{gerald}
\end{document}